\documentclass[review]{siamart0516}

\usepackage{algorithm}
\usepackage{algpseudocode}

\usepackage[T1]{fontenc}
\usepackage[american]{babel}

\usepackage{amsmath}
\usepackage{mathtools}
\usepackage{amssymb}
\usepackage{amstext}
\usepackage{amsfonts}

\usepackage{graphicx}
\ifpdf
  \DeclareGraphicsExtensions{.eps,.pdf,.png,.jpg}
\else
  \DeclareGraphicsExtensions{.eps}
\fi

\newcommand{\ve}[1]{\mathbf{#1}}           
\newcommand{\sv}[1]{\boldsymbol{#1}}   
\newcommand{\m}[1]{\mathbf{#1}}               
\newcommand{\sm}[1]{\boldsymbol{#1}}   
\newcommand{\tr}[1]{{#1}^{\mkern-1.5mu\mathsf{T}}}              
\newcommand{\norm}[1]{||{#1}||}              
\newcommand*{\mvec}{\operatorname{vec}}
\newcommand*{\trace}{\operatorname{trace}}
\newcommand*{\rank}{\operatorname{rank}}
\newcommand*{\diag}{\operatorname{diag}}

\newcommand{\widebar}[1]{\overline{#1}}  

\newcommand{\Had}{\circ}
\DeclareMathOperator*{\lmin}{Minimize}
\DeclareMathOperator*{\argmin}{arg\,min}
\DeclareMathOperator*{\argmax}{arg\,max}

\newcommand*{\intersect}{\cap}
\newcommand*{\union}{\cup}

\let\emptyset\varnothing

\newcommand{\field}[1]{\mathbb{#1}}
\newcommand{\Reals}{\field{R}}

\makeatother
\begin{document}
\title{Euclidean distance matrix completion and point configurations from the minimal spanning tree}
\author{Adam Rahman \and R. Wayne Oldford
}
\headers{MINIMAL SPANNING TREE COMPLETIONS \& CONFIGURATIONS}{Adam Rahman \& R. Wayne Oldford}

\maketitle

\begin{abstract}
The paper introduces a special case of the Euclidean distance matrix completion problem (edmcp) of interest in statistical data analysis where only the minimal spanning tree distances are given and the matrix completion must preserve the minimal spanning tree.   Two solutions are proposed, one an adaptation of a more general method based on a dissimilarity parameterized formulation, the other an entirely novel method which constructs the point configuration directly through a guided random search.  These methods as well as three standard edcmp methods are described and compared experimentally on real and synthetic data.  It is found that the  constructive method given by the guided random search algorithm clearly outperforms all others considered here.  Notably, standard methods including the adaptation force peculiar, and generally unwanted, geometric structure on the point configurations their completions produce.
\end{abstract}



\section{Introduction} \label{sec:introduction}
In applications as varied as network localization, ultrasound tomography, molecular configuration, and statistical data analysis, the objective is to determine a set of locations of points in $p$ dimensions, a point configuration, given only a subset of the distances between points.   This problem has come to be known as the Euclidean distance matrix completion problem and has typically been formulated as an optimization problem, particularly as a semidefinite programming problem -- e.g. see \cite{johnson1995connections,krislock2012euclidean}.  

In this paper, the problem of interest is the special case where the only distances not missing are those which lie on the minimal spanning tree of the point configuration.  The objective becomes not just to produce a point configuration which completes the Euclidean distance matrix but one which also preserves the minimal spanning tree.

Although much information is lost on the geometric configuration of the points when only these distances are available, much is still preserved.  For example, the minimal spanning tree is the backbone of the classic hierarchical clustering method  ``single-linkage'' (e.g. see \cite{gower1969minimum}) as well as more recent improvements such as ``DBSCAN'' \cite{Ester96} or ``runt-pruning'' \cite{Stuetzle03}.  The minimal spanning tree is also used to establish much of the geometry in a scatterplot of points that can be interest to data analysts such as whether the data form ``stringy'', ``sparse'', or ``clumpy'' configurations as defined  by  \cite{wilkinson2005graph}.  
 
A common statistical data analysis application is to take high-dimensional data and embed it in a more manageable lower dimensional space - dimensionality reduction.  In this case all distances and even all point locations could be known for a space of dimension $d >> p$ but one wants a representation of the point configuration in the smaller $p$ dimensional space.  
An early very familiar example would be principal component analysis due to \cite{hotelling1933analysis} where the intent is to preserve in the smaller dimensional space most of the data variation seen in the higher dimensional space.
Similarly, completion of a Euclidean distance matrix containing only the minimal spanning tree distances can reduce the dimensionality of the data while preserving the minimal spanning tree and hence its own important statistical structure. 

In other applications, such as  methodological investigations, it can be of interest to generate synthetic point configurations in $p$ dimensions that have some specified minimal spanning tree.

Euclidean distance matrix completion from a minimal spanning tree is an unusual completion problem.  In most completion problems, the proportion of distances missing from the matrix is relatively small and other features of graphs such as cliques (e.g. see \cite{krislock2010explicit}) can be exploited to solve the problem.  In contrast, completion from any spanning tree has nearly all of the distances missing (e.g. for an $n \times n$ matrix the proportion missing is about 89\% for $n=10$,   99\% for $n = 100$, and 99.9\% for $n=1000$).  It also contains little to no graphical structure to exploit.  

In this paper, we apply some general Euclidean distance matrix completion methods, namely those of Alfakih et. al. \cite{Wolkowicz}, Trosset \cite{Trosset}, and Fang \& O'Leary \cite{OLeary}, to the minimal spanning tree completion problem.  
We also present two completion solutions developed for this problem.  One is a tailoring of the method of \cite{Trosset}, the other is a constructive solution which uses a guided random search to produce the point configuration (and hence the completion).   

All methods are compared experimentally in Section \ref{sec:expt} where two different point configurations are used: the first is real, namely the well-known 4-dimensional Anderson Iris data as given in \cite{Fisher1936use}; the second is synthetic, being several sets of data generated uniformly within a hypercube of varying dimension.  In all cases, the minimal spanning tree is produced and all methods are assessed on performance and on quality of the completions and point configurations they produce.  Additionally, the three general methods are compared in terms of complexity and reconstruction accuracy when the proportion of distances missing varies from 5\% to 95\%.

Somewhat surprisingly, the constructive method outperforms all other methods by every measure especially on computational speed.   Of particular interest are the geometries of the point configurations produced by the various methods.  The three general methods introduce extraneous geometric structure which is generally undesirable in a statistical data analysis.  The minimal spanning tree method adapted from \cite{Trosset} does as well, though not nearly as badly as the other methods.  In marked contrast, the constructive methods introduces very little extraneous geometry in the reconstructed point configuration from its completions. 

The paper is organized as follows.  Section \ref{sec:edm} introduces the notation and provides a brief overview of the Euclidean distance matrix completion problem and the connection to point configurations which are the ultimate goal.  There also, the three completion methods are described and related.  Section \ref{sec:mst} then introduces the problem of minimal spanning tree completion as a subset of the general problem.  The two new methods for solving this problem are introduced here and the algorithms used described  and formally presented.  The experimental results are given in Section \ref{sec:expt};  the Iris data are treated first in Section \ref{sec:Iris} and the uniform data in Section \ref{sec:uniform}.  Concluding remarks are given in Section \ref{sec:Conclusions}.

\section{Euclidean distance matrices and point configurations}\label{sec:edm}
Imagine that we have $n$ points $\ve{x}_1, \ldots , \ve{x}_n \in \Reals^p$ centred so that $\sum_{i=1}^n \ve{x}_i = \ve{0}$ and we denote by $\m{X} = \tr{[\ve{x}_1, \ldots, \ve{x}_n]}$ the $n \times p$ real matrix whose $i$th row is $\tr{\ve{x}}_i$.    Then $\m{D} = [d_{ij}]$ is the $n \times n$ matrix of  {\em squared} Euclidean distances 
\begin{equation}
\label{eq:dij}
d_{ij} = \norm{\ve{x}_i - \ve{x}_j }^2 = \tr{\ve{x}}_i \ve{x}_i +   \tr{\ve{x}}_j \ve{x}_j  - 2  \tr{\ve{x}}_i \ve{x}_j =  g_{ii} + g_{jj} - 2 g_{ij}
\end{equation}
where  $\m{G} =[ g_{ij}] = \tr{\m{XX}}$ is the $n \times n$ {\em Gram} matrix.   The above relationship can be used to define two operators relating the matrices.
Following the notation of \cite{dokmanic2015euclidean} the first is the Euclidean distance matrix operator 
\begin{equation}
\label{eq:edm}
edm(\m{X}) = \ve{1} ~\tr{diag(\m{X}\tr{\m{X}})} - 2 \m{X} \tr{\m{X}} + diag(\m{X}\tr{\m{X}}) ~\tr{\ve{1}}
\end{equation}
which takes the matrix of points $\m{X}$ and turns it into a Euclidean distance matrix.  Here $diag({\m{A}})$ is the vector whose elements are the diagonal entries of the square matrix $\m{A}$ and $\ve{1}$ is a vector of ones.  For convenience, the relationship is also written in terms of an operator on the Gram matrix
\begin{equation}
\label{eq:calK}
{\cal K}(\m{G}) 
= \ve{1} \tr{\ve{g}} - 2 \m{G} + \ve{g} \tr{\ve{1}} 
\end{equation}
where $\ve{g} = diag(\m{G}) = \tr{(g_{11}, g_{22}, \ldots, g_{nn})}$.
The $edm$ operator of  (\ref{eq:edm}) allows us to find a Euclidean distance matrix from any set of points $\ve{x}_1, \ldots, \ve{x}_n$ or equivalently using (\ref{eq:calK}) from the Gram matrix $\m{G}$.

The reverse is also true, namely that we can determine a point configuration $\ve{x}_1, \ldots, \ve{x}_n  \in \Reals^p$ from the Gram matrix $\m{G}$ which in turn can be determined from the Euclidean distance matrix $\m{D}$.  

To see this we first note that equations (\ref{eq:edm}) and (\ref{eq:calK}) imply that we can write the Gram matrix as 
\begin{equation}
\label{eq:gram}
\m{G} = \frac{1}{2} \left( \ve{1} \tr{\ve{g}} -  \m{D} + \ve{g} \tr{\ve{1}} \right).
\end{equation}
If $\ve{g}$ is a function only of $\m{D}$, then an eigen decomposition of the symmetric $\m{G} = \m{U}\sv{\Lambda}\tr{\m{U}}$ can be used to find a point configuration $\m{X} = \m{U}\sv{\Lambda}^{\frac{1}{2}}$ where \\ $\sv{\Lambda}^{\frac{1}{2}} = diag(\sqrt{\lambda_1}, \ldots, \sqrt{\lambda_p}, 0, \ldots,0)$ is a diagonal matrix containing the non-zero eigen-values $\lambda_1 \ge \cdots \ge \lambda_p > 0$ of $\m{G}$.  

To show $\ve{g}$ is a function only of $\m{D}$ we note that
 the $\ve{x}_i$s are centred,  and so we can average over the second index of equation (\ref{eq:dij})  to get
\[ 
\widebar{d}_{i+} =  \frac{1}{n}\sum_{j=1}^n d_{ij} =  \tr{\ve{x}}_i \ve{x}_i +  \frac{1}{n}\sum_{j=1}^n  \tr{\ve{x}}_j \ve{x}_j 
\]
and
\[
\widebar{d}_{++} =  \frac{1}{n}\sum_{i=1}^n\widebar{d}_{i+} = 2 \times  \frac{1}{n}\sum_{i=1}^n\tr{\ve{x}}_i \ve{x}_i .
\]
Together these equations 
yield 
\[ 
g_{ii} = \norm{\ve{x}_i}^2 = \widebar{d}_{i+} - \frac{1}{2} ~ \widebar{d}_{++} 
\] 
which shows that the vector $\ve{g} $ is  a function only of $\m{D}$.  
The relation is more usefully represented in matrix form as
\begin{equation}
\label{eq:g}
\begin{array}{rcl}
\ve{g} \tr{\ve{1}}& =  &\frac{1}{n}\m{D}\ve{1}\tr{\ve{1}} - \frac{1}{2}  \ve{1} \widebar{d}_{++}\tr{\ve{1}} \\
&&\\
& =  &\frac{1}{n}\m{D}\ve{1}\tr{\ve{1}} - \frac{1}{2}     \ve{1} \frac{1}{n}\tr{\ve{1}} \m{D}\ve{1}\frac{1}{n}\tr{\ve{1}} \\
&&\\
& =  &\m{D}\m{H}- \frac{1}{2}  \m{H} \m{D}\m{H} \\
\end{array}
\end{equation} 
where  $\m{H} = \frac{1}{n}\ve{1}\tr{\ve{1}} = \ve{1}(\tr{\ve{1}}\ve{1})^{-1}\tr{\ve{1}}$ is an orthogonal projection matrix.

Equations  (\ref{eq:gram}) and (\ref{eq:g}) together show that the Gram matrix can be directly determined from the Euclidean distance matrix as
\begin{equation}
\label{eq:G}
\m{G} = -\frac{1}{2} \left(\m{I} - \m{H}\right) \m{D} \left(\m{I} - \m{H}\right)= -\frac{1}{2} \m{P} \m{D}  \m{P}.
\end{equation}
The second equality is simply a notational change to highlight the orthogonal projection matrix $\m{P} = \m{I} - \m{H}$ which centres vectors to have zero average entries (applied to both row and column vectors in equation \ref{eq:G}).   An eigen decomposition of the Gram matrix $\m{G}$  will give the point configuration $\m{X}$  known as classical multidimensional scaling (or principal coordinates).

Just as ${\cal K}(\m{G})$ of equation (\ref{eq:calK}) maps a Gram matrix to a Euclidean distance matrix, the  linear operator 
\begin{equation}
\label{eq:calT}
{\cal T}(\m{D} ) =  -\frac{1}{2} \m{P} \m{D}  \m{P}
\end{equation}
maps a Euclidean distance matrix to a Gram matrix.  Note that $\m{G} = \m{X}\tr{\m{X}} = {\cal T}(\m{D})$ is a positive semidefinite matrix.  Moreover, since $\m{P} \ve{1} = \ve{0}$, it follows that $\tr{\ve{z}}\m{D}\ve{z} \le 0$ for any $\ve{z} \in \Reals^n$ satisfying $\tr{\ve{z}}\ve{1}=0$.

More formally, let ${\cal S}_n$ denote the set of real $n \times n$ symmetric matrices and of these let 
\[{\cal S}_n^+ = \{ \m{S} \in {\cal S}_n : \m{S} \succeq 0 \}\]
denote those which are positive semidefinite (denoted $\m{S} \succeq 0$).  

Similarly let ${\cal G}_n  = \{ \m{G} \in  {\cal S}_n :   \m{G}\ve{1}=0 \}$ and ${\cal D}_n = \{ \m{D} \in  {\cal S}_n :   \diag(\m{D}) = \ve{0} \}$ having the following subsets of interest:
\begin{align*}
{\cal G}_n^+ & = \{ \m{G} \in  {\cal G}_n :  \m{G} \succeq 0 \}  ~\mbox{ and }\\
{\cal D}_n^-& = \{ \m{D} \in  {\cal D}_n :   \tr{\ve{z}}\m{D}\ve{z} \le 0 ~ \mbox{ whenever }~ \tr{\ve{z}}\ve{1}=0\} .
\end{align*}
Each of the above subsets is a positive cone in ${\cal S}_n$ of some intrinsic interest:  ${\cal D}_n$ is the set of hollow matrices,  ${\cal D}_n^-$ the set of Euclidean distance matrices, ${\cal G}_n$ is the set of symmetric centred matrices, and ${\cal G}_n^+ $ the set of Gram matrices.

Relations between these sets provide a means for moving back and forth between different representations of the same problem. 
\cite{johnson1995connections}, \cite{OLeary} show that  ${\cal G}_n^+={\cal T}({\cal D}_n^-)$, ${\cal D}_n^-= {\cal K}({\cal G}_n^+)$ and that ${\cal T}$ on ${\cal D}_n^-$ is the inverse of ${\cal K}$ on ${\cal G}_n^+$. 
These lead to several ways in which to characterize a Euclidean distance matrix (EDM).  For example, a hollow matrix $\m{D} \in {\cal D}_n$ is an EDM if and only if $\m{G}={\cal T}(\m{D}) \in {\cal G}_n^+$ (i.e. $\m{G}$ is positive definite).  Alternatively, $\m{D} \in {\cal D}_n$ is an EDM if and only if there exists $\m{G} \in  {\cal G}_n^+$ for which $\m{D} = {\cal K}(\m{G})$.  


\subsection{The completion problem}\label{sec:completion}
Suppose that $\m{D}$ has missing entries and we would like to recover a point configuration in dimension $p$ whose distances match the known distances.   
 
To be concrete, let $\m{A} = [a_{ij}]$ be a matrix where $a_{ij}$ is $1$ if  the squared distance $d_{ij}$ is known and zero otherwise.  The matrix $\m{A}$ can be thought of as an adjacency matrix defining a graph connecting the points in the configuration whenever the distance between them is known.  It is assumed that this graph is connected, otherwise the points separate into disconnected groups with no information to locate the relationship between groups.

This matrix can also be used as a mask to select the elements of a Euclidean distance matrix through a Hadamard (elementwise) product $\m{A} \Had \m{D}$ whose only non-zero values will be the known values of $\m{D}$.   This allows the completion problem to be expressed as the solution to an optimization problem such as  
\begin{equation}
\label{eq:emdcp}
\sm{\Delta}_0 = \argmin_{\sm{\Delta} \in {\cal D}_n^-}~ \norm{\m{A} \Had (\m{D} - \sm{\Delta})}^2_F 
\end{equation}
where $\norm{\m{M}}_F = \norm{\mvec(\m{M})} = \sqrt{\trace(\tr{\m{M}}\m{M})} $ denotes the Frobenius norm of a matrix $\m{M}$.   

With a solution $\sm{\Delta}_0$ in hand, we need only find the Gram matrix $\m{G} = {\cal T}(\sm{\Delta}_0)$ and then its eigen decomposition will provide the matrix of point configurations $\m{X}$.  Of course, since $\m{G} = \m{X} \tr{\m{X}}$,  an equivalent point configuration will be $\m{X}\m{O}$ for any orthogonal $p \times p$ matrix $\m{O}$ where $p$ is the {\em embedding dimension} of $\sm{\Delta}_0$. (Note that the embedding dimension is equal to the rank of $\m{G}$ and so is at most $n-1$.)

Note that the problem could be recast in terms of finding the Gram matrix as
\begin{equation}
\label{eq:edmcpG}
\m{G}_0 = \argmin_{\m{S} \in {\cal G}_n^+}~ \norm{\m{A} \Had (\m{D} - {\cal K}(\m{S}))}^2_F 
= \argmin_{\m{S} \in {\cal G}_n^+}~ \norm{\m{A} \Had  {\cal K}(\m{G} -\m{S})}^2_F 
\end{equation}
where $\m{G} = {\cal T}(\m{D})$, the second of these following from the linearity of the operators.

Various solutions to the problem have been proposed in the literature, each one deriving from a different specification of the constrained minimization. 
In what follows, we consider three such solutions.

\subsubsection{Semi-definite programming formulation}\label{sec:sdp}

Note that ${\cal T}$ of equation (\ref{eq:calT})  pre- and post-multiplies the matrix $\m{G}$  by the orthogonal projection matrix $\m{P}$ of rank $n-1$.  Writing $\m{P} = \m{V} \tr{\m{V}}$ where $\m{V}$ is an $n \times (n-1)$ matrix satisfying $\tr{\m{V}}\sv{1} = \sv{0}$ and $\tr{\m{V}}\m{V} = \m{I}_{n-1}$,  \cite{Wolkowicz} introduce the composite operators
\[ {\cal T}_{V} (\m{D}) = \tr{\m{V}}{\cal T}(\m{D})\m{V}  = -\frac{1}{2}\tr{\m{V}}\m{D}\m{V}
\]  
and
\[ {\cal K}_{V} (\m{S}) = {\cal K}(\m{V}\m{S}\tr{\m{V}})
\]  
and show that ${\cal S}_{n-1}= {\cal T}_{V}({\cal D}_n)$, ${\cal D}_n={\cal K}_{V}({\cal S}_{n-1})$,  and that ${\cal T}_{V}$ and ${\cal K}_{V}$ are inverses of each other on these sets.

With these operators equation (\ref{eq:edmcpG}) can be recast as
\begin{equation}
\label{eq:comp}
\m{S}_0 = \argmin_{\m{S} \in  {\cal S}_{n-1}^{+}}  \norm{\m{A}  \Had  (\m{D} - {\cal K}_{V}(\m{S}))}_{F}^{2} 
\end{equation}
which gives the solutions $\sm{\Delta}_0 = {\cal K}_{V}(\m{S}_0)$, $\m{G} = \m{V}{\cal T}_{V}(\m{D}) \tr{\m{V}}$, and point configuration $\m{X}$.

The search space ${\cal S}_{n-1}^+$ can however be broadened.  \cite{Wolkowicz} show that the composite operators also relate the set ${\cal D}_n^-$ of Euclidean distance matrices to the less restrictive set ${\cal P}_{n-1}^+$ of  real $(n-1) \times (n-1)$ positive semidefinite matrices $\m{P}$.  This allows the problem to now be cast as one of semidefinite programming, namely
\begin{align*}
\lmin &  ~~f(\m{S}) = \norm{\m{A}  \Had  (\m{D} - {\cal K}_{V}(\m{S}))}_{F}^{2}\\
\mbox{subject to:} & ~~~ a(\m{S}) = \ve{c} \\
                              & ~~~\m{S} \succeq 0.
\end{align*}
The elements of $\m{S}$ are any real values and the vector valued function $a(\cdot)$ and constant vector $\ve{c}$ can be used to enforce the constraints on the known squared distances, e.g. ${\cal K}_{V}(\m{S})_{ij} = d_{ij}$.

A major advantage of this semidefinite programming formulation is that it is a convex optimization problem, meaning that any local solution is also a global solution - a property that neither of the other algorithms we consider possesses. However, as pointed out in \cite{OLeary},  this convexity comes at a cost.  The number of parameters is $\mathcal{O}(n^{2})$ which becomes very large for even moderately sized problems.

This formulation also does not allow specification of the embedding dimension $p$, so that the resulting rank of $\m{G}$ and hence $\m{X}$ could be large. A numerical rank approach is proposed by \cite{Wolkowicz} in which very small eigen values are discarded to determine the point configuration with little change to the value of the objective function.

\texttt{Matlab} code is provided by \cite{Wolkowicz} for their algorithm, which we denote by \texttt{SDP} to emphasize its formulation as the solution to a semidefinite programming problem.

\subsubsection{Non-convex position formulation}\label{sec:npf}
Oftentimes we wish to specify the embedding dimension $p$ for the point configuration $\m{X}$.  
The Gram matrix formulation of equation \eqref{eq:edmcpG} can be adapted to a specific embedding dimension $p$ as
\begin{align}
\label{eq:gpcf}
\lmin_{\m{G} \in {\cal G}_n^+} 
&  ~~  \norm{\m{A}  \Had  (\m{D} - {\cal K}(\m{G}))}_{F}^{2}\\
\mbox{subject to:} & ~~~\rank(\m{G}) = p. \nonumber
\end{align}
Fixing the embedding dimension allows the problem to be written as
\begin{equation}
\label{eq:pcf}
\lmin_{\m{X} \in {\Reals}^{n \times p}}   ~~f_{A, D}(\m{X}) = 
         \norm{\m{A}  \Had  (\m{D} - {\cal K}(\m{X}\tr{\m{X}}))}_{F}^{2}.
\end{equation}
The search space here is much smaller (only $np$ values in $\m{P}$) than with the semidefinite programming formulation ($\mathcal{O}(n^{2})$).
Moreover, it is unrestricted.   Although the Gram matrix $\m{G}= \m{X}\tr{\m{X}} \succeq 0$, it need not be in ${\cal G}^+_n$ since there is no longer any restriction that  $\m{G}\ve{1} = \ve{0} = \tr{\m{X}}\ve{1}$.
Adding any constant vector $\ve{c}$ to all points $\ve{x}_i$ in $\m{X}$ makes no change to $\m{D}$ nor to the value of  ${\cal K}(\m{X}\tr{\m{X}})$ and so no change to the objective function $f_{A, D}(\m{X})$.

The objective function of (\ref{eq:pcf}) can now be written in terms of the unconstrained point vectors $\ve{x}_i$ as
\begin{equation}
\label{eq:pcf:objectivefn}
f_{A, D}(\m{X}) = \sum_{i=1}^n \sum_{j=1}^n a_{ij}^2 (d_{ij} - \norm{\ve{x}_i - \ve{x}_j}^2)^2
\end{equation} 
which may now be minimized using more standard numerical optimization procedures.  Unfortunately this formulation is no longer convex and any minimization algorithm may converge to a local solution rather than the global solution.

An algorithm that increases the chances of arriving at a global maximum is introduced by \cite{OLeary}.  The algorithm depends on several insights into the nature of the completion problem.  First, a candidate distance matrix may always be generated as ${\cal T}(\m{B})$ for any hollow matrix $\m{B} \in {\cal D}_n$.  Wherever the squared distances $d_{ij}$ of $\m{D}$ are known, they should be matched by those entries $b_{ij}$ of $\m{B}$.  For any $(i,j)$ for which $d_{ij}$ is not known $a_{ij}$ is zero and the graph given by taking $\m{A}$ to be an adjacency matrix can be traversed to find the distance of the shortest path between $i$ and $j$.  This distance provides a conservative choice for $b_{ij}$ when $a_{ij}=0$.  

This distance can be large and so is made less conservative by considering the number of edges along the shortest path. Test cases of Euclidean distance matrix completion problems are investigated in \cite{OLeary} to arrive at a more useful choice for these $a_{ij}$.  In the end, they choose to randomly generate values for these $b_{ij}$ as 
\[b_{ij} = \frac{f_{ij}}{s_{ij}} \]
where $f_{ij}$ is the shortest path distance between $i$ and $j$, and $s_{ij}$ is a realization from a $N(1.5, 0.009)$ distribution truncated to be between $0$ and the number of segments, $k_{ij}$, in the shortest path. A large number of such matrices $\m{B}$ are randomly generated in \cite{OLeary}, choosing to start the optimization using whichever ${\cal T}(\m{B})$ produced the smallest value of $f_{A, D}(\m{X})$.

Since the global minimum will yield $f_{A, D}(\m{X})=0$ it is easy to tell when the optimization has settled on a local minimum.  Two means of getting off the local solution are offered in \cite{OLeary}.  First, they allow ``stretching'' whereby all distances are increased by a common scale and the optimization started again.  Secondly, they allow the embedding dimension to be increased so that a global optimization can be more easily found and then restart the optimization from this solution after it has been reduced to the desired embedding space (e.g. by projection).

These random starts and means of expanding the problem to avoid local minimum are used in conjunction with a careful combination of optimization procedures to arrive at a solution. We denote the algorithm of \cite{OLeary} by \texttt{NPF} to emphasize it as a solution to the non-convex position formulation problem.  

\subsubsection{Dissimilarity parameterized formulation}\label{sec:dbf}
The problem is formulated in \cite{Trosset} by focussing first on $n \times n$ dissimilarity matrices
\[{\cal C}_n = \{ \sm{\Delta} \in {\cal D}_n :   \sm{\Delta} = [ \delta_{ij}] ~\text{ has } \delta_{ij} \ge 0 \}\]
which clearly contains the Euclidean distance matrices  ${\cal D}_n^- \subset {\cal C}_n \subset {\cal D}_n$.  This set is then restricted to the subset 
\begin{equation}
\label{eq:Cn}
{\cal C}_n ( \m{A} \Had \sm{\Delta}^{\star} ) = \{ \sm{\Delta} \in {\cal C}_n :   \m{A} \Had \sm{\Delta} =  \m{A} \Had  \sm{\Delta}^{\star}  \}
\end{equation}
containing all completions of the partial dissimilarity matrix $\m{A} \Had  \sm{\Delta}^{\star}$.  Let ${\cal D}_n^-(p) \subset {\cal D}_n^-$ denote the $n \times n$ Euclidean distance matrices from $p$ dimensional point configurations.
Then whether there exists a solution to the completion problem having embedding dimension $p$ is equivalent to whether ${\cal C}_n ( \m{A} \Had \m{D} ) \intersect {\cal D}_n^-(p) \ne \emptyset$.  
This intersection is shown to be non-empty in \cite{Trosset} if and only if the following minimization 
\begin{align}
\label{eq:dis1}
\lmin_{\m{G}, \sm{\Delta}}
&  ~~  \norm{\m{G} - {\cal T}(\sm{\Delta}))}_{F}^{2}\\
\mbox{subject to:} & ~~~\m{G} \in {\cal S}_n^+ ~\text{ and }~  \rank(\m{G}) \le p \nonumber \\
                              & ~~~\sm{\Delta} \in {\cal C}_n ( \m{A} \Had \m{D} ) \nonumber 
\end{align}
has a global minimum of zero.  This formulation is similar to that of equation (\ref{eq:gpcf}).   Note however, that the masking matrix $\m{A}$ and the known distances $\m{A} \Had \m{D} $ now appear as part of the constraints given by the set of allowed dissimilarity matrices ${\cal C}_n ( \m{A} \Had \m{D} )$.

Letting
\[
{\cal C}_n ( \sm{\Delta}^L, \sm{\Delta}^U) = \{ \sm{\Delta} \in {\cal C}_n :   \delta_{ij}^L \le \delta_{ij} \le \delta_{ij}^U   \}
\]
where $\sm{\Delta}^L = [ \delta_{ij}^L]$ and $\sm{\Delta}^U = [ \delta_{ij}^U]$ are both in ${\cal C}_n $, we see that formulation (\ref{eq:dis1}) is a special case of 
\begin{align}
\label{eq:dis2}
\lmin_{\m{G}, \sm{\Delta}}
&  ~~  \norm{\m{G} - {\cal T}(\sm{\Delta}))}_{F}^{2}\\
\mbox{subject to:} & ~~~\m{G} \in {\cal S}_n^+ ~\text{ and }~  \rank(\m{G}) \le p \nonumber \\
                              & ~~~\sm{\Delta} \in {\cal C}_n ( \sm{\Delta}^L, \sm{\Delta}^U) \nonumber 
\end{align}
where $ \sm{\Delta}^L$ and $ \sm{\Delta}^U$ are fixed matrices chosen such that $ \m{A} \Had \sm{\Delta}^L   = \m{A} \Had  \sm{\Delta}^U$ and whenever $a_{ij} = 0$ then $\delta_{ij}^L =0$ and $\delta_{ij}^U =+\infty $.  A huge advantage of this formulation is that lower and upper bounds can be added to further restrict the space of possible solutions -- e.g. using the structure of the graph given by $\m{A}$ upper bounds $\delta_{ij}^U$  can be determined for all unknown $\delta_{ij}$ by simply invoking the triangle inequality.

By minimizing first over $\m{G}$, \cite{Trosset97} shows the problem (\ref{eq:dis2}) can be reduced to
\begin{align}
\label{eq:dpf}
\lmin_{\sm{\Delta}}
&  ~~  F_p({\cal T}(\sm{\Delta}))\\
\mbox{subject to:} &  ~~~\sm{\Delta} \in {\cal C}_n ( \sm{\Delta}^L, \sm{\Delta}^U) \nonumber 
\end{align}
when 
\[
F_{p}(\m{S}) = \sum_{i=1}^{p}(\lambda_{i} - \lambda_{max})^{2} + \sum_{i=p+1}^{n}\lambda_{i}^{2}
\]
for any  $\m{S} \in {\cal S}_n$ whose eigen-values are $\lambda_1 \ge \cdots \ge \lambda_n$ and so $\lambda_{max} = \lambda_1$.  

An efficient algorithm is provided in \cite{Trosset} for solving problem \eqref{eq:dpf} which we will denote by \texttt{DPF} to emphasize its dissimilarity parameterized formulation.

\section{Completions from a minimum spanning tree}\label{sec:mst}

Our interest lies in completions of a Euclidean distance matrix when the graph given by $\m{A}$ determines a spanning tree on its $n$ nodes and the known distances given by $\m{A} \Had \m{D}$ are also the minimum spanning tree distances of the completion.  

All methods in the previous section can be applied when only the first condition holds, though there may be computational challenges given that  $1 - \frac{2}{n} $ of the distances are missing.  None are \emph{specifically} designed to ensure that the minimal spanning tree remains unchanged in the completion.  The goal here is to find {\em mst-preserving} completions.

More formally, let ${\cal A}_n \subset {\cal D}_n $ denote the set of symmetric $n \times n$ adjacency matrices,  ${\cal A}_n^{\star} \subset {\cal A}_n $ the subset corresponding to spanning trees, and  ${\cal A}_n^{\star} (\m{A}) \subset {\cal A}_n^{\star}$ the set of spanning tree adjacency matrices restricted to a subgraph of $A$ (i.e. $\m{A} - \m{A}^{\star} \in {\cal A}_n$ whenever  $\m{A}^{\star} \in {\cal A}_n^{\star}$).  
Further let $amst(\m{A}, \sm{\Delta}) = \{\m{A}_1, \ldots, \m{A}_k\}$ denote the set of adjacency matrices $\m{A}_i \in {\cal A}_n^{\star} (\m{A})$ which produce a minimal spanning tree for the graph given by $\m{A}$ and $\sm{\Delta} \in {\cal C}_n$.  This set will typically be a singleton, but could be larger whenever there are tied values within a dissimilarity matrix $\sm{\Delta}$.  The product $\m{A}_i \Had \sm{\Delta}$ will determine a minimal spanning tree for any $\m{A}_i \in amst(\m{A}, \sm{\Delta})$.

Now let
\[
{\cal M}_n ( \m{A}, \m{A}^{\star}, \sm{\Delta}^{\star} ) = \{ \sm{\Delta} \in {\cal C}_n : amst(\m{A}, \sm{\Delta}) =   amst(\m{A}^{\star}, \sm{\Delta}^{\star}) \}
\]
denote all those dissimilarity matrices $ \sm{\Delta}$ which with a  given  adjacency matrix $\m{A}$ will have the same minimal spanning tree adjacency matrix set as that for the target $\sm{\Delta}^{\star}$ and $\m{A}^{\star}$.  When non-empty, the set ${\cal M}_n ( \m{A}, \m{A}^{\star}, \sm{\Delta}^{\star} )  \subset {\cal C}_n$ is a positive cone (proof in Appendix).
The set of mst-preserving completions of $\m{A}^{\star} \Had \sm{\Delta}^{\star}$ is now simply the intersection of ${\cal C}_n(\m{A}^{\star} \Had \sm{\Delta}^{\star})$ and ${\cal M}_n ( \m{A}, \m{A}^{\star}, \sm{\Delta}^{\star} )$ \emph{when} $\m{A} = \m{K} = \ve{1}\tr{\ve{1}} - \diag(\ve{1}) $ is the adjacency matrix for a complete graph on $n$ nodes.   

Analogous to Equation \eqref{eq:Cn}, we define this set to be
\[ 
{\cal M}_n ( \m{A} \Had \sm{\Delta}^{\star} ) = \{ \sm{\Delta} \in {\cal C}_n(\m{A} \Had \sm{\Delta}^{\star}):   amst(\m{K}, \sm{\Delta}) =  amst(\m{A},  \sm{\Delta}^{\star})  \}
\]
and note now that mst-preserving completions with embedding dimension $p$ exist if and only if ${\cal M}_n ( \m{A} \Had \m{D} ) \intersect {\cal D}_n^-(p) \ne \emptyset$.  Note also that the minimal spanning tree fixes only very few dissimilarities/distances (viz. $n-1$) and leaves a great many to be determined (viz. $(n-1)(n-2)/2$).  Moreover, the fixed distances are the smallest that produce a spanning tree.  The set ${\cal M}_n ( \m{A} \Had \sm{\Delta}^{\star})$ is very large.

As in \cite{Trosset}, an mst-preserving completion problem in $p$ dimensions can now be expressed as
\begin{align}
\label{eq:mstcp}
\lmin_{\m{G}, \sm{\Delta}}
&  ~~  \norm{\m{G} - {\cal T}(\sm{\Delta}))}_{F}^{2}\\
\mbox{subject to:} & ~~~\m{G} \in {\cal S}_n^+ ~\text{ and }~  \rank(\m{G}) \le p \nonumber \\
                              & ~~~\sm{\Delta} \in {\cal M}_n ( \m{A} \Had \m{D} ) \nonumber 
\end{align}
achieving a zero global minimum.  This differs from the formulation \eqref{eq:dis1} only in restricting $\sm{\Delta}$ to ${\cal M}_n ( \m{A} \Had \m{D} )$, a  subset of ${\cal C}_n ( \m{A} \Had \m{D} )$, which suggests that the methods of \cite{Trosset} could be adapted to find an mst-preserving completion by solving a minimization problem.    That is,  we might reduce \eqref{eq:mstcp} to \eqref{eq:dpf} exactly as before, provided we ensure that $\sm{\Delta}^L$ and $\sm{\Delta}^U$ are chosen so as to ensure the other constraints hold.

\subsection{Judicious choice of bounds}\label{sec:dbflb}
In \cite{Trosset} solving \eqref{eq:dpf} only specified $ \delta_{ij}^L =0$; no greater lower bound is used.  If lower bounds can be determined so that the minimal spanning tree is maintained then solving  \eqref{eq:dpf}  will also solve \eqref{eq:mstcp}.  This observation suggests that an adaptation of the \texttt{DPF} algorithm in \cite{Trosset} using the correct non-zero lower bounds will produce an mst-preserving completion.

Algorithm \ref{alg:mstLB} 
\begin{algorithm}[htbp]  
    \caption{MST lower bounds algorithm}
     \label{alg:mstLB}
     {\footnotesize 
    \algsetblock[Name]{Setup}{EndSetup}{1}{0.5cm}               
    \algsetblock[Name]{Structures}{EndStructures}{3}{0.5cm}               
    \begin{algorithmic}[0] 
        \Structures
            \State tree: $T =(V, E)$ is a spanning tree with vertex set $V =\{v\}$ and edge set $E=\{e\}$; 
            \State edges:  $e$ will be a set of two indices $\{i, j\} = nodes(e)$ and have a weight $wt(e) = \delta_{ij} \ge 0$;
            \State $\sm{\Delta}^L = [\delta_{ij}^L]$ is the matrix of dissimilarity lower bounds to be determined;\\
            
        \Procedure{splitTree}{$T, splitEdge$}\Comment{$T$ is a spanning tree}
             \State $restEdges  \gets edges(T) - \{splitEdge\}$\Comment{Remove $splitEdge$ from the edge set of $T$}
             \State $ (v_1, v_2) \gets nodes(splitEdge)$
             \State $ V_1 \gets \{v_1\}$;  $E_1 \gets \{e \in restEdges : v_1 \in nodes(e) \}$\Comment{$E_1$ could be empty}
             \State $ V_2 \gets \{v_2\}$;  $E_2 \gets \{e \in restEdges : v_2 \in nodes(e) \}$\Comment{$E_2$ could be empty}

             \While{$restEdges \ne \emptyset$}
                 \State {$e \gets restEdges[1]$}
                 \State $restEdges \gets restEdges - \{e\}$
                 \If {$nodes(e) \intersect nodes(E_1)$}
                      \State $E_1 \gets E_1 \union \{e\}$
                 \Else 
                      \State $E_2 \gets E_2 \union \{e\}$
                 \EndIf
             \EndWhile
             \State $V_1 \gets V_1 \union nodes(E_1)$; $V_2 \gets  V_2 \union nodes(E_2)$
             \State {\bf return} $\{T_1 \coloneqq  (V_1 , E_1 ), ~~ T_2 \coloneqq  (V_2 , E_2 )\}$ \Comment{Return the two trees}
        \EndProcedure \\
        
         \Procedure{MSTLowerBounds}{$T, \sm{\Delta}^L$}\Comment{Recursively determines the lower bounds}
            \If{$edges(T) \ne \emptyset$}\Comment{Ensure there are edges left in $T$}
                 \State  $ maxEdge \gets \argmax_{e \in edges(T)} wt(e)$\Comment{Split on the biggest edge}
                 \State  $ Trees  \gets \textsc{splitTree}(T, maxEdge)$
                 \For{$v_1 \in nodes(Trees[T_1])$}
                     \For{$v_2 \in nodes(Trees[T_2])$}
                         \State $\sm{\Delta}^L[v_1,v_2]  \gets \sm{\Delta}^L[v_2,v_1] \gets wt(maxEdge)$\Comment{Set the lower bound}
                     \EndFor
                 \EndFor
                \For{$Tree \in Trees$}
                   \State $\sm{\Delta}^L \gets$ \textsc{MSTLowerBounds}$(Tree, \sm{\Delta}^L)$ \Comment{Recursively get lower bounds}
                \EndFor 
            \EndIf
            \State {\bf return} $\sm{\Delta}^L$ \Comment{Return the matrix of lower bounds}
        \EndProcedure
    \end{algorithmic}
    } 
\end{algorithm}
is used in conjunction with the \texttt{DPF} algorithm to construct lower bounds for all distances that will preserve the minimum spanning tree in the completion.  The insight in constructing the lower bound is drawn from single-linkage clustering.  Every edge in a spanning tree separates the vertices into two different groups, depending on which points remain connected to either one vertex or the other of that edge.  Because the tree is a minimum spanning tree, if we select the largest edge, then the distance between any vertex of one group and any vertex of the other group must be at least as large as that of the the largest edge.  This gives a lower bound for these distances that will preserve that edge in the minimum spanning tree.  The same reasoning is applied recursively to each separate group, thus producing a lower bound on all edges.

With these lower bounds computed in advance, we need only use these in Trosset's algorithm to find a solution to  \eqref{eq:mstcp}.  We denote the method as  \texttt{DPFLB} to indicate that it is simply
\texttt{DPF} with these lower bounds.

\subsection{A constructive solution}\label{sec:con}
For any $\m{A}$ and $\m{D}$,  the set $\sm{\Delta} \in {\cal M}_n ( \m{A} \Had \m{D} )$ is large; for a connected graph it is at its largest when $\m{A}$ specifies a minimal spanning tree.  In this case, it should be possible to find a completion in $\sm{\Delta} \in {\cal M}_n ( \m{A} \Had \m{D} )$.  Here, we simply construct such a completion by locating points $\ve{x}_i  \in  \Reals^p$ ($p$ being the embedding dimension) one at a time,  while checking that the (partial) minimal spanning tree is preserved as each point is added.

More formally, we define $\m{X}_k = \tr{[\ve{x}_1, \ldots, \ve{x}_k]}$ to be a $k \times p$ matrix whose rows are point locations $\ve{x}_1, \ldots, \ve{x}_k \in \Reals^p$.  The locations are chosen so that the minimal spanning tree from the Euclidean distances of these $k$ locations in $\Reals^p$ is identical to that of $k$ connected nodes from the minimal spanning tree $\m{A} \Had \m{D}$.  
The matrices $\m{X}_1, \m{X}_2, \ldots, \m{X}_n$ are constricted by growing (and preserving) the minimal spanning tree one node at a time.  The distance matrix from $\m{X}_n$ provides an mst-preserving completion.

The construction begins by choosing the node of maximal degree from the minimal spanning tree of $\m{A} \Had \m{D}$ and locating it at $\ve{x}_1 = \ve{0}$.   The second node to locate will be that of maximal degree amongst those connected to the first.  If the dissimilarity between these two nodes is, say, $\delta_{12}$, then the location of $\ve{x}_2$ is chosen at random from a uniform distribution on the surface of a sphere $S^{p-1}$  in $\Reals^p$ of radius $\sqrt{\delta_{12}}$ centred at $\ve{x}_1$ (assuming squared distances for the completion matrix).  The two points $\ve{x}_1$ and $\ve{x}_2$ are trivially a subtree of the minimal spanning tree.   The remaining nodes with connections to $\ve{x}_1$ are then added in similar fashion.  

As each location is proposed, its (squared) distance to all other placed points is calculated and the resulting distance matrix checked to see whether the minimal spanning tree (so far) is preserved.  If it is, the point is accepted; if not, points are generated until one is acceptable.   When new nodes are added, they are chosen amongst those without locations that share an edge in the minimal spanning tree with nodes already located.  At each step, available nodes of highest degrees are added before nodes of low degree in  $\m{A} \Had \m{D}$; since these are harder to place, they appear earlier.   
Algorithm \ref{alg:constructive} describes the method in detail. 
\begin{algorithm}[htbp]
    \caption{Constructive Completion Algorithm}
    \label{alg:constructive}
    {\footnotesize
    \algsetblock[Name]{Setup}{EndSetup}{1}{0.5cm}               
    \algsetblock[Name]{Structures}{EndStructures}{7}{0.5cm}               
    \begin{algorithmic}[0] 
        \Structures
            \State $T_V =(V, E)$ is a tree spanning the vertex set $V = \{1, \ldots, n\}$ with edge set $E$; 
            \State $\m{A} = [a_{ij}]$ and $\m{A}^{\ast} = [a_{ij}^{\ast}]$  are $n \times n$ symmetric adjacency matrices;
            \State $\sm{\Delta} = [\delta_{ij}]$, is an  $n \times n$ symmetric squared dissimilarity matrix;
            \State $\m{X}$ is an $n \times p$ point configuration matrix to be constructed.
  \\          
            \Procedure{mstConfigure}{$\m{A}, \sm{\Delta},~  \texttt{maxIn = 100},~  \texttt{maxOut = 100}$}  
                 \State  $\texttt{nTries} \gets 0$; ~~$\texttt{Converged?}\gets \texttt{FALSE}$;  $T_V \gets (V,E) \gets tree(\m{A})$;
           
             \Repeat 
                 \State $\texttt{nTries} \gets \texttt{nTries} + 1$
                 \State $\m{X} \gets \m{0}$; $\sm{A} \gets \m{0}$;   \Comment{Initialization}
                  \State $i \gets \argmax_{j \in V} degree(node(j)) $ \Comment{Start at any maximal degree node $i$ in $T_V$}
                  \State  $P \gets  \{i\}$  \Comment{Initial vertex set}
                  \State $T_{P} \gets (P, \emptyset)$\Comment{Root the tree}
                  \State $\texttt{Grow?} \gets \texttt{TRUE}$ \Comment{Keep growing flag}
             \\
             \While{$\texttt{Grow?} $} 
                  \State $(g, B) \gets \textsc{getBuds}(T_{P}, T_V) $ \Comment{grow $B \subset{V}$, $B \intersect P = \emptyset$ from $g \in T_P$}
        
                    \For{$b \in B$}
                \State $(T_P,  \sm{\Delta}, \m{A}, \m{X}, \texttt{Converged?}, \texttt{Grow?}) \gets \textsc{growTree}(g, b, T_P, T_V, \sm{\Delta}, \m{A}, \m{X}, \texttt{maxIn}$) 
                   \EndFor
             \EndWhile         
             \Until $\texttt{nTries} > \texttt{maxOut}$ or $\texttt{Converged? == TRUE}$   
                
          \State {\bf{return $( \m{X}, \texttt{Converged?})$}}\Comment{Return the point configuration}
        \EndProcedure \\
        
        \Procedure{growTree}{$i, j, T_P, T_V, \sm{\Delta}, \m{A}, \m{X}, \texttt{maxTries} = 100$} 
            \State $\texttt{Placed?} \gets \texttt{Converged?} \gets \texttt{FALSE}; ~~\texttt{nTries} \gets 0;\m{A}^{\ast} \gets \m{A};~~ \sm{\Delta}^{\ast} \gets \sm{\Delta}; $\Comment{Initialization}  
            \Repeat 
              \State $\texttt{nTries} \gets  \texttt{nTries} +1 $
              \State $\ve{z} \sim \texttt{Uniform}(S^{p-1})$\Comment{Generate a random direction vector}
              \State $\ve{x}_{j} \gets \ve{x}_{i} + \ve{z} \times \sqrt{\delta_{ij}} $\Comment{Propose the point}
                  \For {$ k \in P$} \Comment{Try values for all placed nodes}
                      \State  $\delta_{jk}^{\ast} \gets\delta_{kj}^{\ast} \gets \norm{\ve{x}_{k} - \ve{x}_{j}}^{2}$
                      \State $a_{jk}^{\ast} \gets a_{kj}^{\ast} \gets 1$
                  \EndFor
                  \If {$\textsc{mst}(\m{A}^{\ast} \Had \sm{\Delta}^{\ast}) \subset T_V$}   $\texttt{Placed?} \gets \texttt{TRUE}$\Comment{Preserves the MST?}
                  \EndIf
            \Until{\texttt{Placed?} or  $\texttt{nTries} > \texttt{maxTries}$} 
                  
            \If{\texttt{Placed?}} \Comment{Accept the point}
            \State $\m{X}[j, ] \gets \tr{\ve{x}_j}; ~~\m{A} \gets \m{A^{\ast}}$;  $~~\sm{\Delta} \gets \sm{\Delta^{\ast}}; $
            \State $nodes(T_P) \gets P \union \{j\}$; $edges(T_P) \gets edges(T_P) \union \{(i,j), (j,i)\}$
                \If{$T_P = T_V$} $\texttt{Converged?} \gets \texttt{TRUE}$
                \EndIf
            \EndIf
            \State {\bf{return} $(T_P,  \sm{\Delta}, \m{A}, \m{X}, \texttt{Converged?}, \texttt{Placed?})$}
        \EndProcedure
        \\
        \Procedure{getBuds}{$T_P, T_V$} 
                    \State  $B \gets V - P$
                    \State  $E_{\ast} = \{ (i,j) : i \in P, j \in B, (i,j) \in edges(T_V) \}$ \Comment{edges in $T_V$ connecting $P$ and $B$}
                    \State $g \gets \argmax_{i \in P}  \#\{ e: e \in E_{\ast} ~~\text{and} ~~ i\in e\}$ \Comment{$g \in P$ having most connections to $B$}
                    \State $B \gets \{b : (g,b) \in E_{\ast} \}$\Comment{reduce $B$ to nodes connected to $g$}
                    \State {\bf{return} $(g, B)$}
        \EndProcedure
        
    \end{algorithmic}
    } 
\end{algorithm}

This is a guided random search algorithm, simply denoted as \texttt{C},  that constructs the configuration point by point.

\section{Experimental results}\label{sec:expt}
In Section \ref{sec:edm} three methods of completing Euclidean distance matrices and generating point configurations were described: \texttt{SDP}, \texttt{NPF}, and \texttt{DPF}.   To these two new methods were added in Section \ref{sec:mst}:  \texttt{DPFLB} and \texttt{C}.
In this section, we experimentally assess these five different completion algorithms in a variety of ways.  

Each experiment begins with a known point configuration $\m{X}$ and its Euclidean distance matrix $\m{D}$.   Elements of $\m{D}$ will be removed and each method will be expected to find a completion $\sm{\Delta} = \widehat{\m{D}}$ and a corresponding point configuration $\widehat{\m{X}}$.   The quality of each method is then based on the nearness of $\widehat{\m{X}}$ and $\widehat{\m{D}}$ to the original $\m{X}$ and $\m{D}$, respectively.  For each, the time taken to arrive at a solution is measured or, equivalently, the number of solutions produced in the same time.

We use two different point configurations.  The first is the well known Anderson Iris data as collected by \cite{Anderson:1935} and recorded in  \cite{Fisher1936use}.   The data consist of four measurements (sepal width and length, petal width and length) on each of 150 flowers, 50 from each of three different Iris species (Versicolor, Virginica, and Setosa).  So $\m{X}$ is a $150 \times 4$ matrix in Section \ref{sec:Iris}. 
The second will be a family of synthetic configurations drawn randomly from uniform distributions on unit hypercubes of varying dimensionality.  That is each row location of $\m{X}$ is  randomly generated as $\ve{x}_i \sim U[0,1]^p$ for $i=1, \ldots, n$.  In the experiments which follow $n=100$, $p \in \{2, \ldots, 10\}$ and  for each $p$ five different point configurations are generated (called ``matrix 1'' to ``matrix 5'' in Section \ref{sec:uniform}). 

\subsection{Reconstructing the Iris data}\label{sec:Iris}
For the three methods from Section \ref{sec:edm} matrix completions can be generated for any pattern of missing dissimilarities provided the corresponding graph spans all points.  To investigate the relative performance of these methods we take $\m{D}$ to be constructed from the distances between flowers in the Iris data.  Distances are removed at random and the percentage of distances removed varied.  Each completion method produces $\widehat{\m{D}}$ for several randomly selected patterns of missing distances for each percentage missing.

\subsubsection{Completions as a function of percentage missing}\label{sec:Iris:percentmissing}
For each completion  the time taken is recorded.  To measure the accuracy of each completion  the relative difference in dissimilarities 
\begin{equation}
\label{eq:rdd}
RDD \, = \, \frac{||\m{D} - \m{\widehat{D}}||_{F}^{2}}{||\m{D}||_{F}^{2}}
\end{equation}
is calculated as well.

\begin{figure}[hbt]
\begin{center}
\begin{tabular}{cc}
\includegraphics[width=0.85\textwidth]{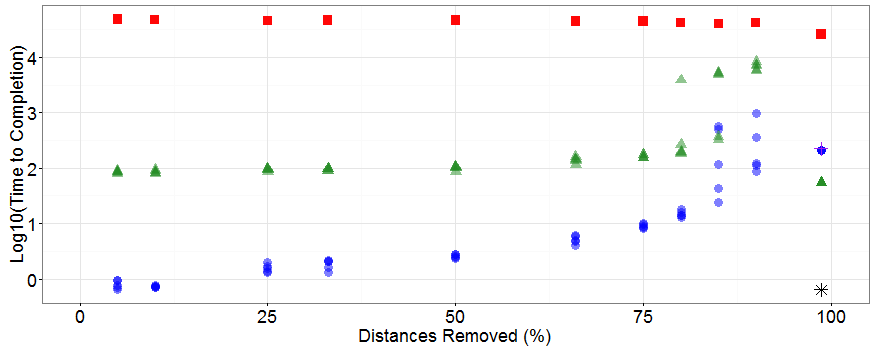} &
\hspace{-0.025\textwidth}
\includegraphics[trim=0cm 0cm 0cm 0cm, clip=TRUE,width=0.12 \textwidth]{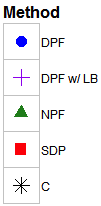} 
\\
\includegraphics[width=0.85\textwidth]{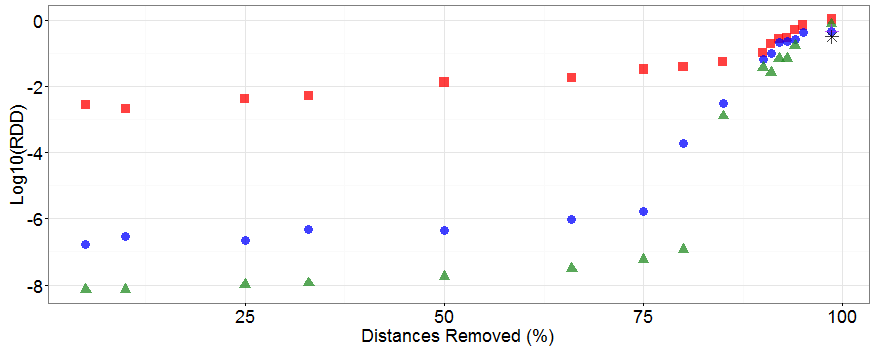} &
\\
\end{tabular}
\end{center}
\caption{\emph{Effect of increasing the percentage of missing distances on each method.  Methods are coded by colour and symbol shape.  For each percentage, several different matrices with different patterns of randomly selected missing data were used.  In the top plot, each point symbol represents the result of one such matrix with one method; alpha blending of colour is used so that places where the values are essentially the same will appear more saturated due to over-plotting and the blending of the colours. In the bottom plot, only the average values are shown.}}
\label{fig:IrisPerformance}
\end{figure}

The top plot of Figure \ref{fig:IrisPerformance} 
shows the effect of percentage missing has on the computational time.  For each percentage, every method but \texttt{SDP} was applied to the same five different incomplete missing at random matrices; the \texttt{SDP} method took so long that it was applied to only the first of the five matrices.  Not surprisingly  computational time increases with the percentage missing.  Comparing methods, we see that the \texttt{DPF} method of  \cite{Trosset} is consistently faster than the \texttt{NPF} method of \cite{OLeary}, which in turn is consistently faster than the \texttt{SDP} method of \cite{Wolkowicz}.  Computational times are given on a logarithmic scale so these differences are substantial.  Some variation in the times taken can also be seen, especially for the larger percentages.  

The minimal spanning tree case has the greatest percentage missing possible and appears at the far right of each plot.    All times to completion here appear to have dropped, with \texttt{NPF} and \texttt{DPF} switching positions to make it fastest of the three methods.  Since this is the minimal spanning tree case, the two methods of Section \ref{sec:mst} can be added.  Not surprisingly, \texttt{DPFLB} which differs from \texttt{DPF} only in having precomputed non-zero lower bounds for every missing distance takes essentially the same time to complete as does \texttt{DPF}.  More interestingly, the constructive method \texttt{C} is orders of magnitude faster than all other methods.

The lower plot of Figure \ref{fig:IrisPerformance} shows the average accuracy with which the various methods reconstruct $\m{D}$.   On this logarithmic scale lower values indicate greater accuracy so the accuracy of every method decreases as the percentage missing increases -- largely because the numerator in \eqref{eq:rdd} has more non-zero entries while the denominator remains unchanged. 
For every percentage up to and including 85\% missing, each method was applied to 50 different missing at random matrices, the exception being \texttt{SDP} which, because of the time required, was only applied to a single matrix each time.  Beyond 85\% only 10 different random matrices were used for each of \texttt{NPF} and \texttt{DPF}.   Again the methods can be ordered: \texttt{NPF} is consistently most accurate and \texttt{SDP} consistently least accurate across all percentages missing.  

For the minimal spanning tree case, there is only one matrix to complete but all methods (with the exception of \texttt{SDP}) have a random step (only if necessary for \texttt{NPF}) and so could produce many potentially different solutions.  Each method was allowed to complete as many completions as possible in the same time taken for \texttt{SDP} to construct one completion -- \texttt{NPF} was able to make 466 completions (all identical since the random step never needed to be invoked), \texttt{DPF} made 124 completions, \texttt{DPFLB} 114, and the constructive method produced a remarkable 39,935 completions!  

\begin{figure}[htbp]
\begin{center}
\includegraphics[width=0.8\textwidth]{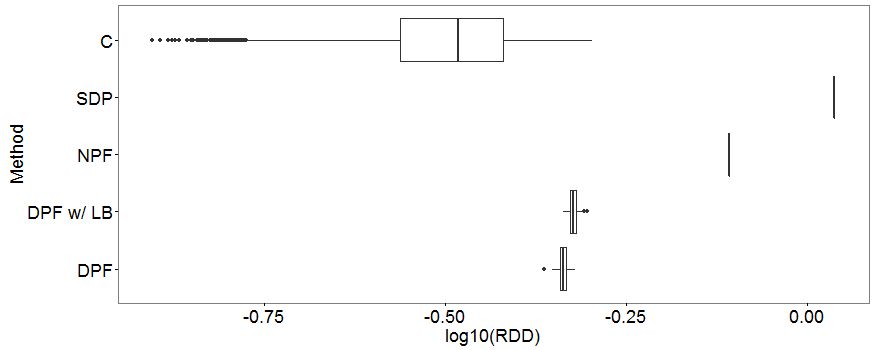}
\end{center}
\caption{\emph{Relative dissimilarity difference when completing the Iris data from its minimal spanning tree distances.  \texttt{C} has 39,935 completions, \texttt{SDP} one,  \texttt{NPF} 466 identical completions, \texttt{DPF}  124 completions, and  \texttt{DPFLB} 114 -- each set was completed in the time taken for a single \texttt{SDP} completion ($> 7$ hours)}}
\label{fig:mstIrisAccuracy}
\end{figure}

Figure \ref{fig:mstIrisAccuracy}
displays boxplots of the $\log_{10} RDD$ for the five methods.  Clearly, the constructive method \texttt{C} nearly always outperforms all of the others.   From most to least accurate there is the constructive method \texttt{C} as most accurate, then \texttt{DPF} closely followed by \texttt{DPFLB}, then  \texttt{NPF}, and finally \texttt{SDP}.  
As was the case with the time to completion, note that again the \texttt{DPF/LB} and \texttt{NPF} switched order in the minimal spanning tree case.

It would seem that \texttt{SDP} performed poorest on both measures, whatever the percentage missing.  
For \texttt{DPF} and \texttt{NPF} one performed better than the other on each measure and neither dominated the other on both.  
When we consider beginning only with the minimal spanning tree distances, the constructive method \texttt{C} performed best on both measures with the improvement in time being  considerable.

\subsubsection{Distances as a function of percentage missing} \label{sec:Iris:Distances}
\begin{figure}[hbtp]
\begin{tabular}{cccl}
\includegraphics[width=0.23\textwidth]{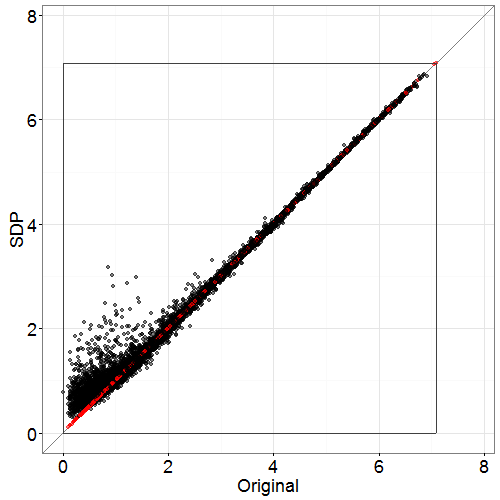}  &
\includegraphics[width=0.23\textwidth]{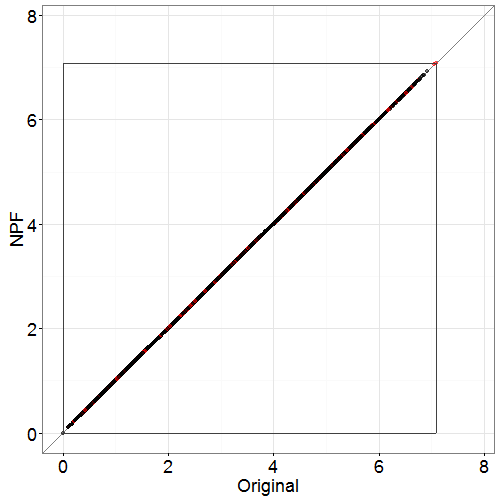}  &
\includegraphics[width=0.23\textwidth]{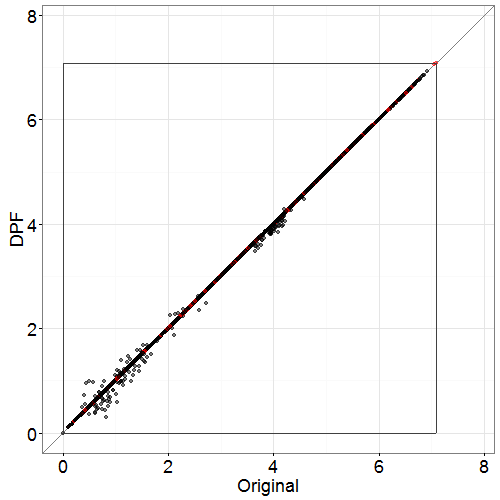} & 85\% removed\\ 
\includegraphics[width=0.23\textwidth]{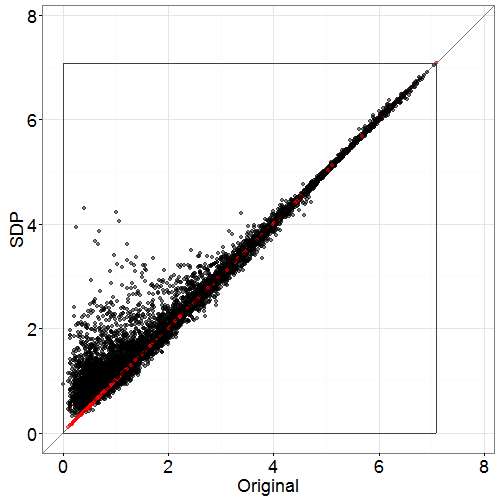}  &
\includegraphics[width=0.23\textwidth]{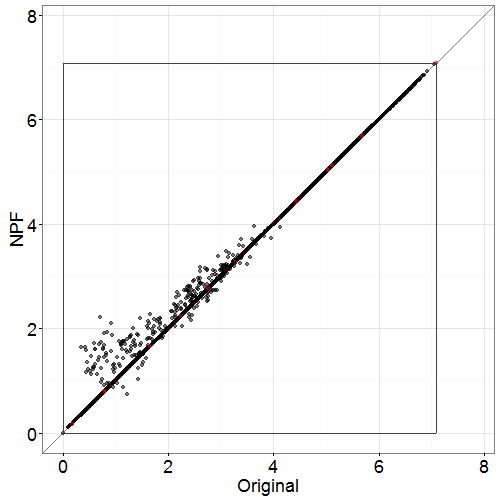}  &
\includegraphics[width=0.23\textwidth]{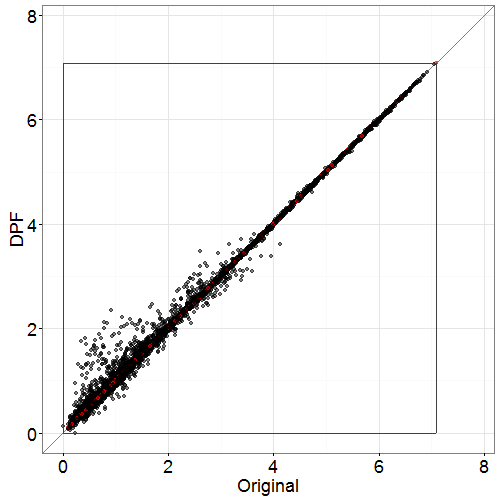} & 90\% removed\\ 
\includegraphics[width=0.23\textwidth]{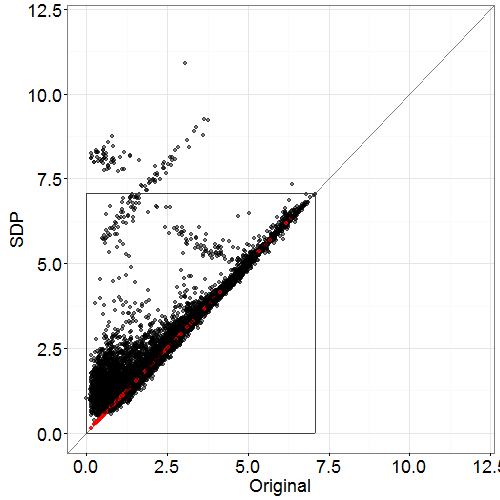}  &
\includegraphics[width=0.23\textwidth]{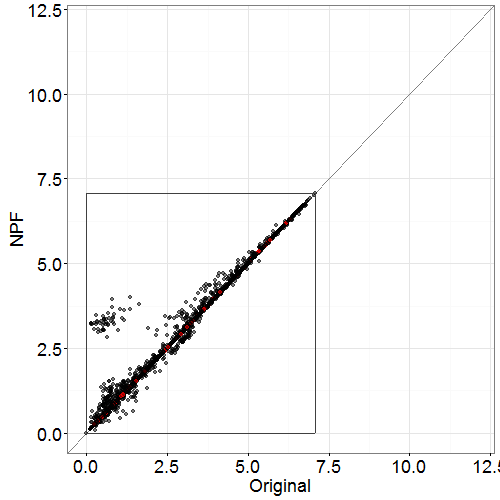}  &
\includegraphics[width=0.23\textwidth]{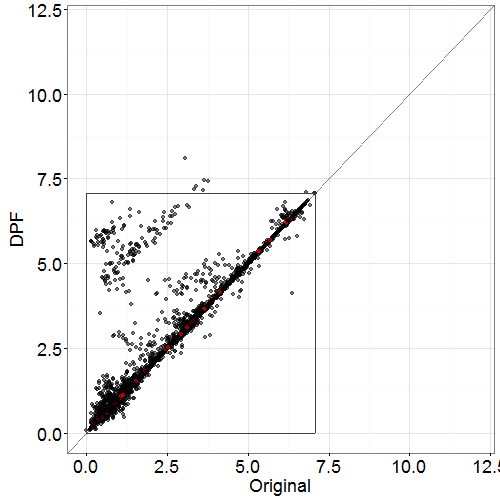} & 92\% removed\\ 
\includegraphics[width=0.23\textwidth]{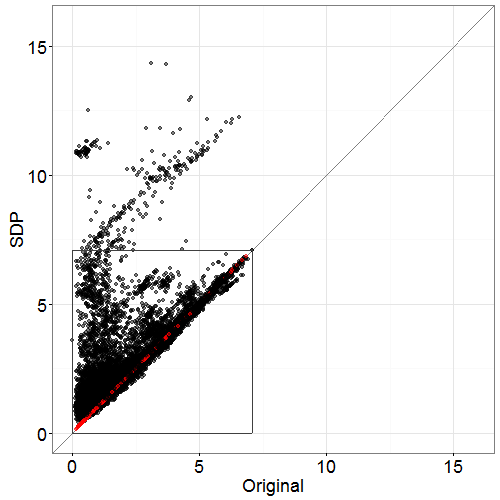}  &
\includegraphics[width=0.23\textwidth]{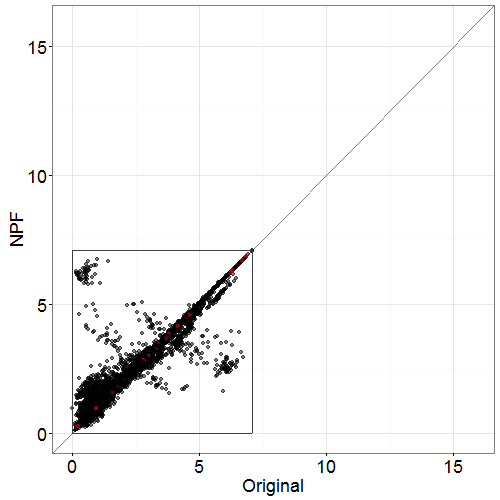}  &
\includegraphics[width=0.23\textwidth]{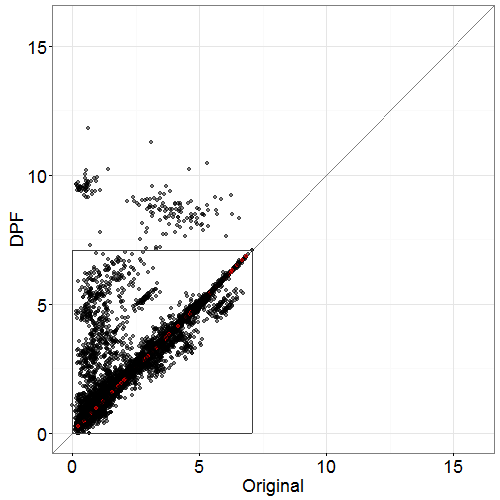} & 94\% removed\\ 
\footnotesize{\texttt{SDP}} &
\footnotesize{\texttt{NPF}} &
\footnotesize{\texttt{DPF}} &\\
\end{tabular}
\caption{\emph{Plots of the pairs of $(d_{ij}, \widehat{d}_{ij})$ for all $i < j$ for a single reconstruction of the Iris distance matrix; red values are the original minimal spanning tree distances.   Perfect reconstruction would be all points on the $y=x$ line;  the box in each plot is the range of the original distances and is identical in size across all plots; scales are identical for the same percentage of missing distances.}
\label{fig:IrisDistancesVaryMissing}}
\end{figure}

Further insight into the three completion algorithms can be had by comparing the reconstructed distances $\widehat{\m{D}}=[\widehat{d}_{ij}]$ with the actual distances $\m{D} = [d_{ij}]$. 
Figure \ref{fig:IrisDistancesVaryMissing}
plots the pairs $(d_{ij}, \widehat{d}_{ij})$ for all $i < j$ for a few of the missing percentages.  Perfect reconstruction would be all points on the $y=x$ line;  the box in each plot is the range of the original distances and is identical in size across all plots; scales are identical for the same percentage of missing distances.  

For each case, \texttt{NPF} outperforms the other two -- all distances appear within the box, appear on either side of the $y=x$ line and is nearer to this line in all cases.   As the percentage missing increases, the reconstruction of all three methods degrades.   Both  \texttt{DPF} and \texttt{SDP} tend to produce ever larger distances in their reconstructions as the percentage increases.  \texttt{DPF} does produce some distances that are smaller than the original too.  In contrast, \texttt{SDP} has a tendency to consistently produce distances that are too large for every percentage,  and produces much larger distances than does \texttt{DPF}.  Overall \texttt{NPF} provides the best reconstructed distances and \texttt{SDP} the worst.

\begin{figure}[hbtp]
\begin{tabular}{ccccc}
\includegraphics[width=0.17\textwidth]{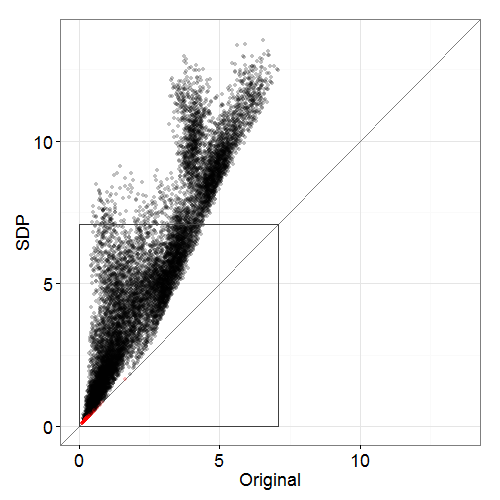}  &
\includegraphics[width=0.17\textwidth]{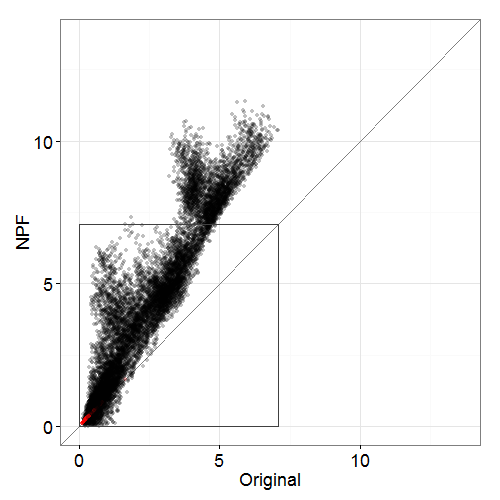}  &
\includegraphics[width=0.17\textwidth]{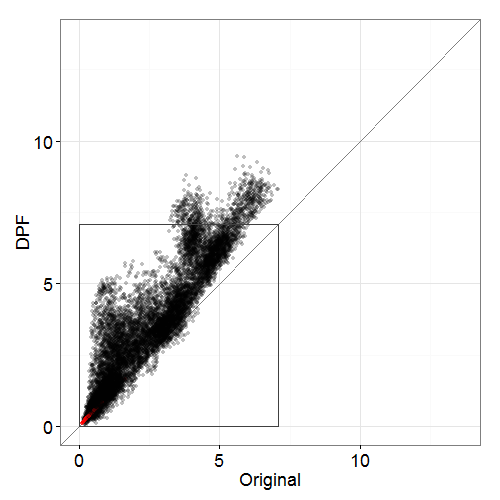} &
\includegraphics[width=0.17\textwidth]{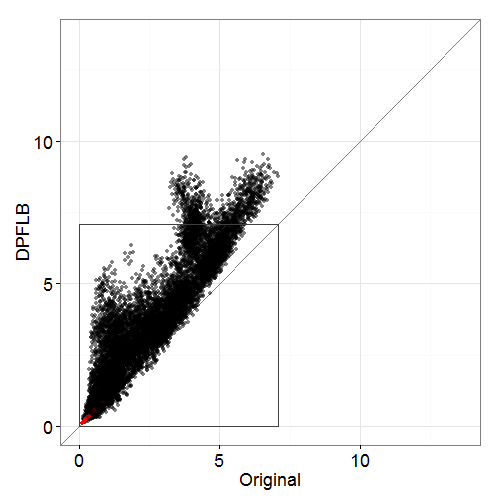} &
\includegraphics[width=0.17\textwidth]{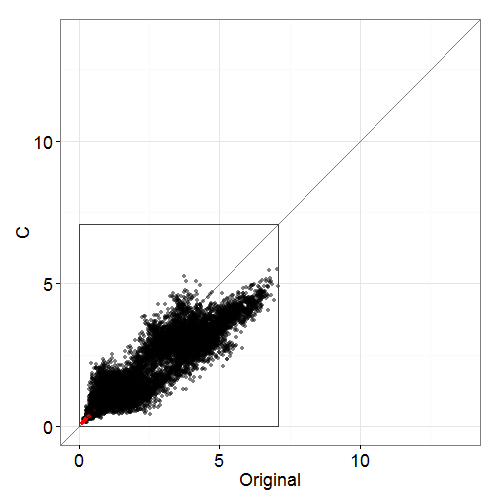} \\ 
\footnotesize{\texttt{SDP}} &
\footnotesize{\texttt{NPF}}&
\footnotesize{\texttt{DPF}} &
\footnotesize{\texttt{DPFLB}} &
\footnotesize{\texttt{C}} 
\end{tabular}
\caption{\emph{Plots of $(d_{ij}, \widehat{d}_{ij})$ for all $i < j$ for a single reconstruction when the matrix to be completed contained only minimal spanning tree distances (shown in red);  the $y=x$ line indicates perfect matching; the box in each plot shows the extent of the original distances.}}
\label{fig:IrisMSTdistances}
\end{figure}

Turning to minimal spanning tree completions, Figure \ref{fig:IrisMSTdistances}
shows the pairs $(d_{ij}, \widehat{d}_{ij})$ for all $i < j$ for a single reconstruction for all five methods.   The shapes of the first four (\texttt{SDP}, \texttt{NPF}, \texttt{DPF}, and \texttt{DPFLB}) are surprisingly similar, each showing three different branches.  All four have almost all distances $\widehat{d}_{ij} > {d}_{ij}$ and many larger than $\max_{ij}{d}_{ij}$, with \texttt{SDP} producing the largest distances, followed by \texttt{NPF}, then  \texttt{DPFLB}.  The last of these has larger distances than those of \texttt{DPF} likely because of the lower bounds which ensure that \texttt{DPFLB} is mst-preserving.  

In marked contrast, the constructive method \texttt{C} produces a completion whose distances $\widehat{d}_{ij}$ are all within the range of the true distances ${d}_{ij}$ and are much more nearly concentrated around the $y=x$ line.  If anything, \texttt{C} seems more inclined to produce smaller distances than are necessary.  This might be corrected by having the vectors generated in Algorithm \ref{alg:constructive} not be  generated on a sphere uniformly but to favour directions that would increase distances to other points already in the tree.

\subsubsection{Reproducing the minimal spanning tree}\label{sec:Iris:repmst}
\begin{figure}[hbt]
\begin{center}
\begin{tabular}{cc}
\includegraphics[height=0.195\textheight]{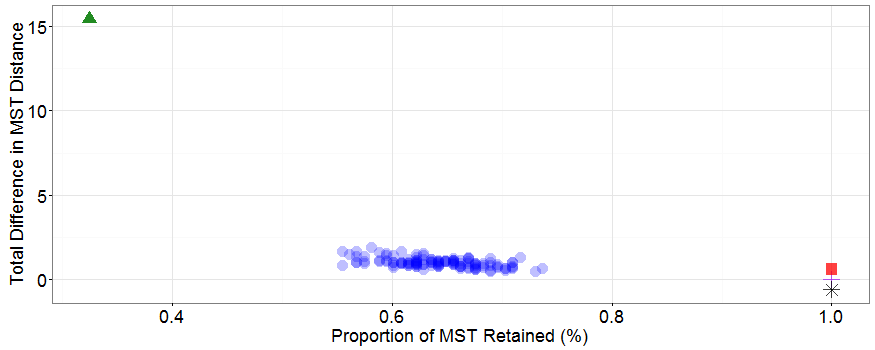} &
\hspace{-0.025\textwidth}\includegraphics[trim=0cm 0cm 0cm 0cm, clip=TRUE,width=0.125 \textwidth]{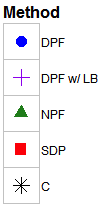} 
\\
\end{tabular}
\end{center}
\caption{\emph{Square of the total difference in minimal spanning tree distances as a multiple of the total squared minimal spanning tree distances versus the proportion of minimal spanning tree edges that were retained.  All completions of all five methods are shown.  The three points at the right have identical values but have been given different vertical positions to better distinguish the points in the plot.}}
\label{fig:IrisMSTPerformance}
\end{figure}

To see how well the various methods reproduced the minimal spanning tree, the spanning tree distances $\widehat{d}_i$ for $i=1, \ldots, n-1$, and adjacency matrix $\widehat{A}$ were determined from each completion $\widehat{\m{D}}$.  Figure \ref{fig:IrisMSTPerformance} 
shows all completions by all methods where each horizontal location is the proportion of edges in $\m{A}$ which also appear in $\widehat{\m{A}}$ and each vertical location is 
\[
\frac{(\sum_{i=1}^{n-1}d_{i} - \sum_{i=1}^{n-1}\widehat{d}_{i})^{2}}{\sum_{i=1}^{n-1}d_{i}^{2}}
\]
where   $d_i$  for $i=1, \ldots , n-1$ are the original minimal spanning tree distances of $\m{A} \Had\m{D}$ and $\widehat{d}_i$ for $i=1, \ldots , n-1$ 
are those from $\widehat{\m{A}} \Had \widehat{\m{D}}$ for that completion.   

The mst-preserving completions are those in the bottom right of Figure \ref{fig:IrisMSTPerformance} -- \texttt{DPFLB},  \texttt{C} , and \texttt{SDP} (all values on both measures but have been separated vertically to better distinguish their shapes from the others methods).   Of these, two were designed to be mst-preserving and so should appear here; each of these two points actually represent hundreds or tens of thousands of completions which must return the same minimal spanning tree.  In contrast, the point for \texttt{SDP} is a singleton point representing the one completion actually constructed for \texttt{SDP}.  This \texttt{SDP} completion has turned out to be mst-preserving, though not by design.  As seen in Figures \ref{fig:IrisDistancesVaryMissing} and \ref{fig:IrisMSTdistances}, \texttt{SDP} tends to produce very large distances in its completions, and these distances increase as the percentage missing increases.  This would explain why \texttt{SDP} is mst-preserving here.  Unfortunately, it also means that given any spanning tree whether minimal or not, this spanning tree will likely be the minimal spanning tree of a completion by \texttt{SDP}.

The 466 identical completions of \texttt{NPF} all appear in the top left corner of Figure \ref{fig:IrisMSTPerformance}  and show \texttt{NPF}  to be the poorest performer in terms of preservation of the minimal spanning tree.  The 124 completions by \texttt{DPF} are spread across the bottom, retaining about 55-75\% of the edges in the minimal spanning tree and matching the distances fairly closely (at least compared to \texttt{NPF}).

\subsubsection{Reproducing the point configurations}\label{sec:Iris:PointConfigs}
We now examine the point configurations produced by all five methods for a single reconstruction from the minimal spanning trees.  

To standardize the comparisons, the Iris data is transformed to its principal directions from a singular value decomposition of $\m{X} = \m{U}\sm{\Lambda}\tr{\m{V}}$ as $\m{X}\m{V}$.  This projects each point $\m{x}_i$ onto a new coordinate system given by the new variates $V1$, $V2$, $V3$, and $V4$. 

\begin{figure}[hbtp]
\begin{tabular}{cc}
\includegraphics[width=0.5\textwidth]{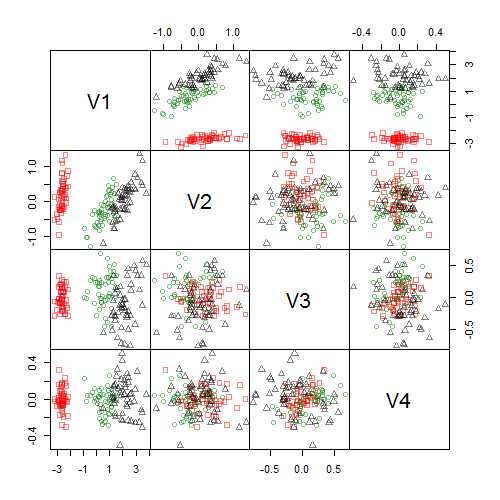} &
\includegraphics[width=0.5\textwidth]{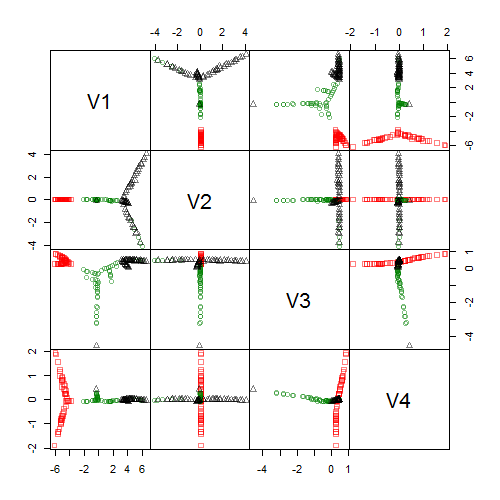} \\
\footnotesize{(a) Original}&
\footnotesize{(b) \texttt{SDP}} \\

\includegraphics[width=0.5\textwidth]{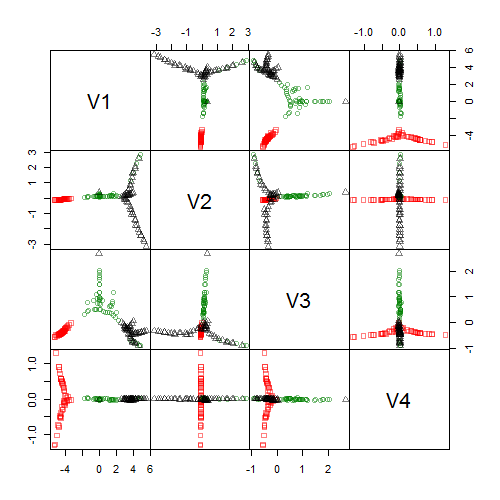} &
\includegraphics[width=0.5\textwidth]{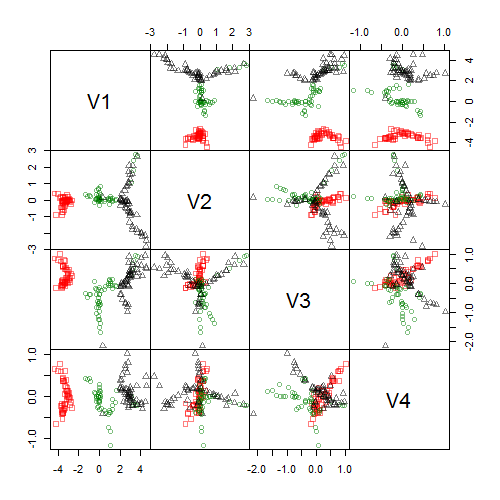} \\
\footnotesize{(c) \texttt{NPF}}&
\footnotesize{(d) \texttt{DPF}} \\

\includegraphics[width=0.5\textwidth]{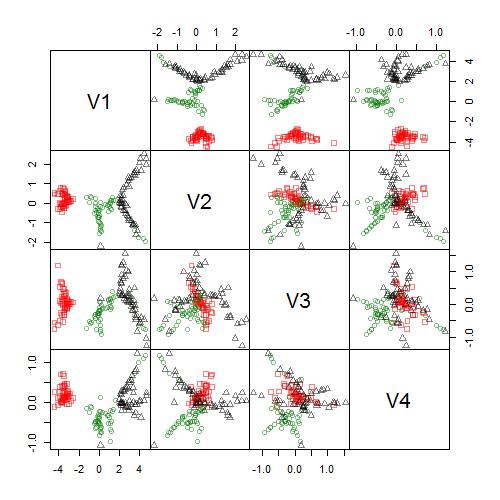} &
\includegraphics[width=0.5\textwidth]{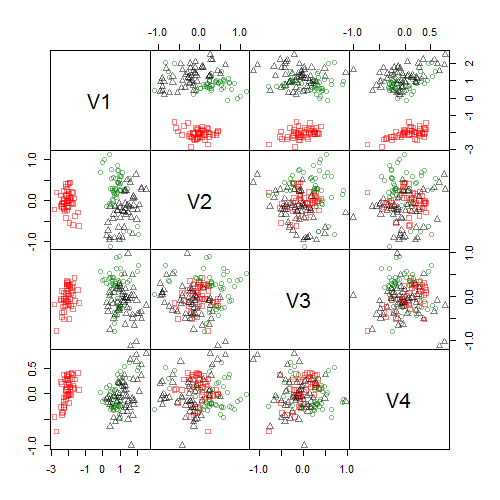} \\
\footnotesize{(e) \texttt{DPFLB}}&
\footnotesize{(f) \texttt{C}} \\
\end{tabular}\caption{\emph{Iris data point configurations reconstructed in the four dimensions given by each configuration's principal coordinates. The three species of flower are distinguished both by colour and by shape of the point symbols.}}
\label{fig:IrisScatters}
\end{figure}

Figure \ref{fig:IrisScatters}(a)
shows the original Iris data in scatterplot matrix for this new coordinate system.  Each scatterplot shows the transformed data for the pair of variates given by the diagonal entries from the same row and column as the scatterplot.  The three different point colours and shapes identify the three different species of Iris.  There are 150 points.

For each completion method, a single completion $\widehat{\m{D}}$ is taken with embedding dimension $p=4$,  the singular value decomposition of its estimated point configuration $\widehat{\m{X}}= \widehat{\m{U}}\widehat{\sm{\Lambda}} \tr{\widehat{\m{V}}}$ determined, and the transformed point configuration $\widehat{\m{X}}\widehat{\m{V}}$ plotted in a scatterplot matrix with the same point symbols.  The transformed point configurations for the five completion methods are shown as Figures \ref{fig:IrisScatters}(b--f).

As can be seen, \texttt{SDP} and \texttt{NPF} produce star-shaped configurations of straight lines.  So too does \texttt{DPF} although the shapes are slightly noisier.   \texttt{DPFLB} is very much like \texttt{DPF} but may be slightly noisier again.  Most striking is the configuration produced by the constructive method \texttt{C}; of the five methods considered, \texttt{C} produces a configuration most like that of the original data.  We note also that the axis ranges in these scatterplots reflect the size of the corresponding singular values from the configurations.  These are consistent with the remarks made earlier about the comparative size of the distances produced by each method in Figure \ref{fig:IrisMSTdistances}.  Again, the configuration of \texttt{C} appears to be closest to the original in this respect as well.  Finally, all completions appear to preserve much of the group separation seen in the original data between the three species of flower.  This is not too surprising since the distances from the minimal spanning tree were given.  The minimal spanning tree is the basis for many hierarchical clustering methods.  The two methods designed to preserve the minimal spanning tree should fare best in maintaining separation of clusters.

\subsection{Reconstructing data $\sim U[0,1]^p$}\label{sec:uniform}
Here we simulate $n$ locations $\ve{x}_i$ independently from $U[0,1]^p$ for varying values of $p$ to give an $n \times p$ point configuration $\m{X}$ within the unit hypercube in $\Reals^p$.
For each $p$, five different matrices of point configurations $\m{X}_1, \ldots, \m{X}_5$ are generated; these will allow us to get a sense of the variability in the results.  

As with the Iris data, the performance of the five methods will be compared but, rather here only completions from the distances of the minimal spanning tree of each configuration are considered.  Of interest then is how this performance might depend on varying dimensionality $p$ rather than the percentage of distances missing.

In what follows, $p  \in \{2,3,4, \ldots, 10\}$ will be used.  The $i$th configuration matrix of dimension $q \le p$ will be constructed as the first $q$ columns of $\m{X}_i$ where $\m{X}_i$ is the $n \times p$ matrix whose rows were independently generated from  a $U[0,1]^{10}$ distribution.  Throughout, we take $n=100$.

\subsubsection{Completions as a function of dimension}\label{sec:Uniform:completions}
\begin{figure}[htbp]
\begin{center}
\begin{tabular}{cc}
\hspace{-0.05\textwidth}
\includegraphics[width=0.9\textwidth]{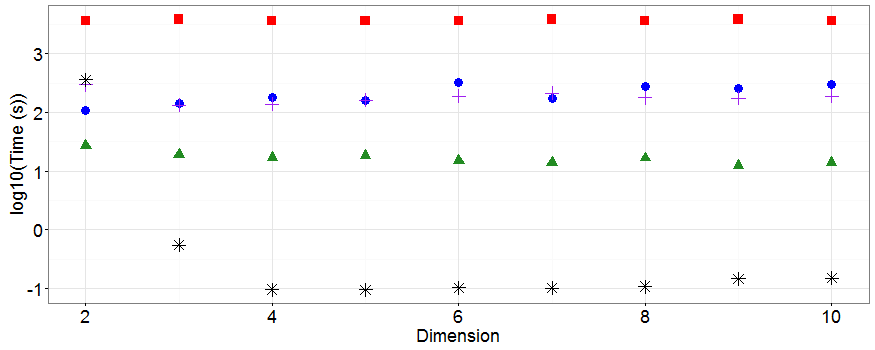} &
\includegraphics[width=0.15\textwidth]{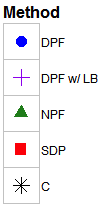}
\\
\hspace{-0.05\textwidth}
\includegraphics[width=0.9\textwidth]{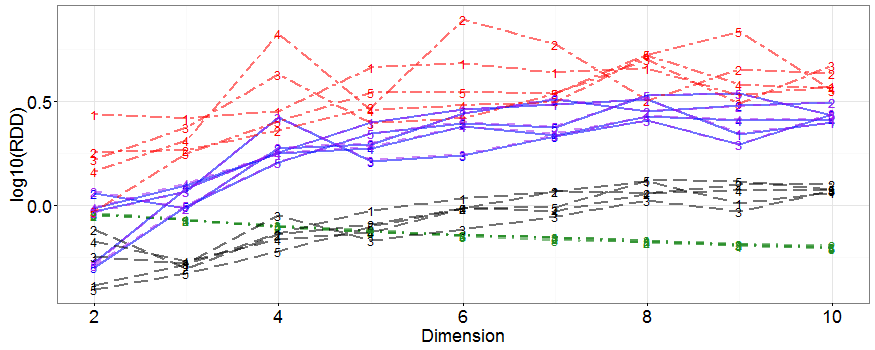} &
\includegraphics[width=0.15\textwidth]{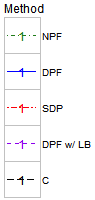}
\end{tabular}
\caption{\emph{Performance of the five completion matrices. Top plot shows base 10 logarithms of the average times to complete ($n=1$ for \texttt{SDP}, $n = 466$ for \texttt{NPF},  $n = 124$ for \texttt{DPF}, $n = 114$ for \texttt{DPFLB}, $n = 39,935$ for \texttt{C}) from a single matrix (``Matrix 1'').  The bottom plot shows  $\log_{10}$ of the average $RDD$s for these completions from all five matrices (marked 1 to 5).}}
\label{fig:UnifPerformance}
\end{center}
\end{figure}

Figure \ref{fig:UnifPerformance}
shows the effect of varying dimension $p$ on the time to completion (top plot) and on the accuracy, as measured by $RDD$ of Equation \eqref{eq:rdd} using all distances, for the three algorithms.   The completion times are from a single matrix (``matrix 1'') for each dimension; the accuracies are shown for each of the five different simulates matrices (marked 1 through 5) with methods distinguished by line type and colour.  Recall that a lower dimensional matrix shares its columns with all higher dimensional ones.  

For time to completion, there is a clear ordering of methods from the least efficient \texttt{SDP} to the several orders of magnitude more efficient \texttt{C}.  Both \texttt{DPF} and \texttt{DPFLB} take about the same time and \texttt{NPF} is second fasted though still two orders of magnitude slower than \texttt{C}.  The only exception to this ordering occurs for $p=2$.  There the two mst-preserving methods are slower than both \texttt{DPF} and \texttt{NPF}; this is likely due to the increased difficulty in finding completions which preserve the minimal spanning tree for uniformly generated data in only two dimensions. 
As $p$ increases, \texttt{SDP} stays relatively constant in computation time, \texttt{NPF} decreases slightly, both \texttt{DPF} and \texttt{DPFLB} tend to increase, and \texttt{C} drops quickly as it becomes easier in larger dimensional spaces to find random directions that work.

In terms of accuracy, all methods degrade as dimensionality increases with the notable exception of \texttt{NPF} whose accuracy improves.   Unfortunately, \texttt{NPF} does not preserve the minimal spanning tree.  The others in order from least to most accurate over all dimensions are \texttt{SDP}, \texttt{DPF} and \texttt{DPFLB} (about the same), and \texttt{C}.  Note again that the logarithm has been taken of the average $RDD$s so these differences can be substantial.

\subsubsection{Distances as a function of dimension}\label{sec:Uniform:Distances}
\begin{figure}[hbtp]
\begin{center}
\begin{tabular}{ccccc}
\includegraphics[width=0.17\textwidth]{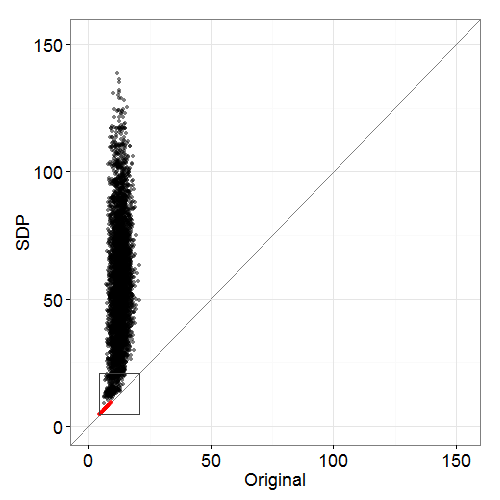} &
\includegraphics[width=0.17\textwidth]{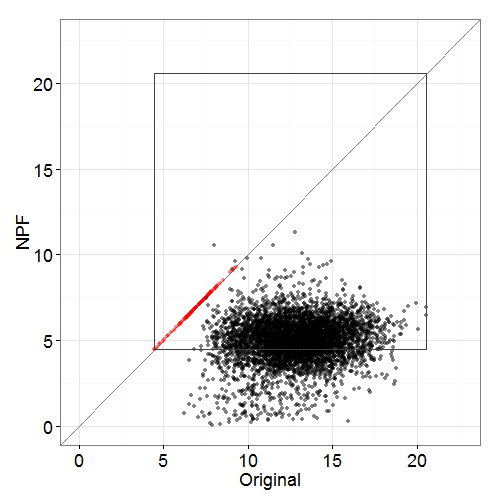} &
\includegraphics[width=0.17\textwidth]{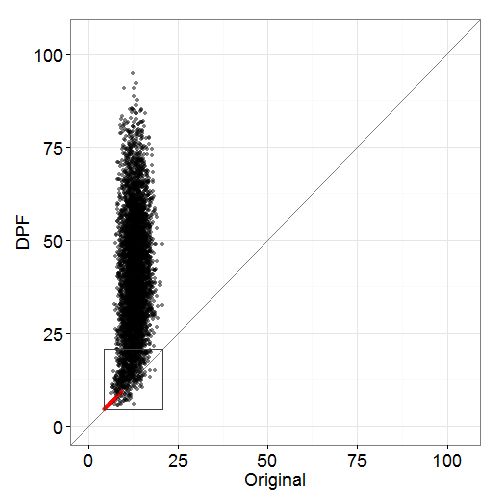} &
\includegraphics[width=0.17\textwidth]{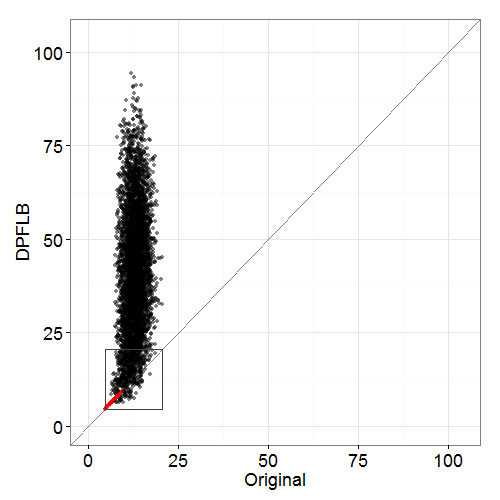}&
\includegraphics[width=0.17\textwidth]{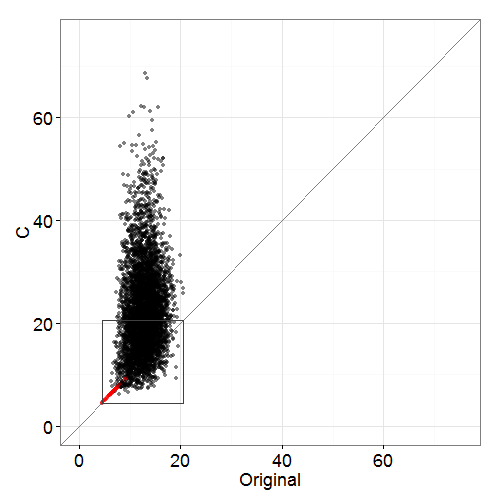} \\ 
\includegraphics[width=0.17\textwidth]{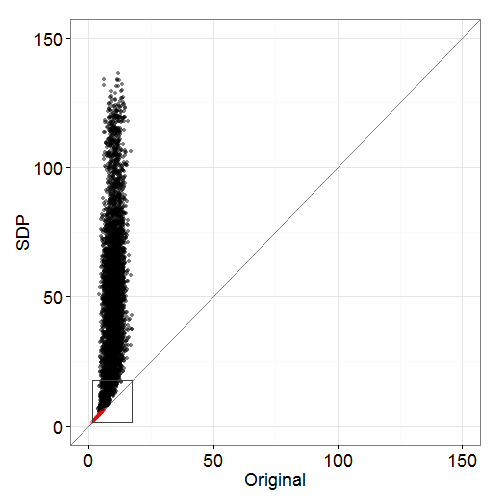} &
\includegraphics[width=0.17\textwidth]{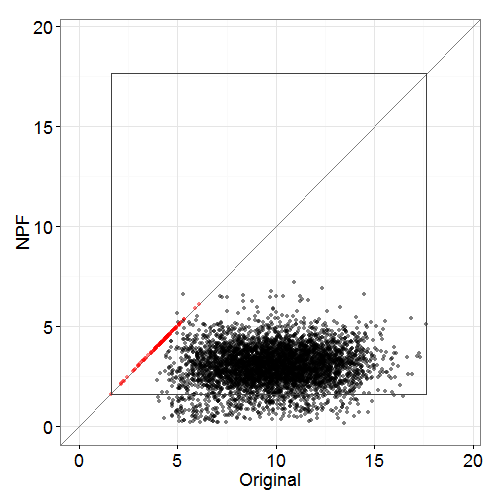}&
\includegraphics[width=0.17\textwidth]{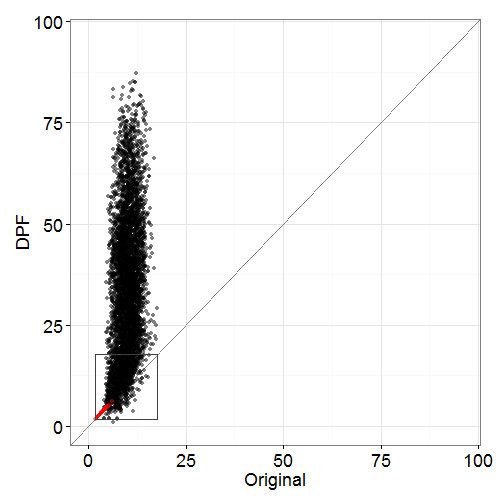}&
\includegraphics[width=0.17\textwidth]{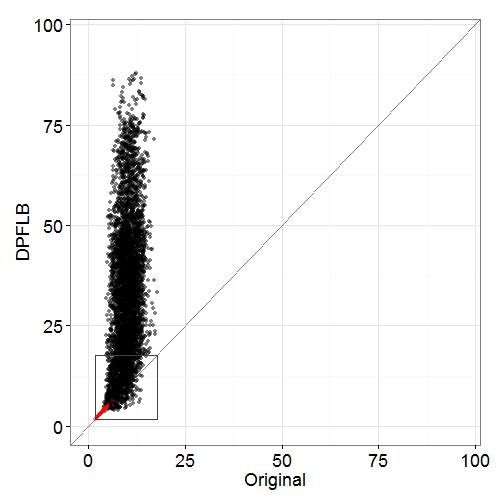}  &
\includegraphics[width=0.17\textwidth]{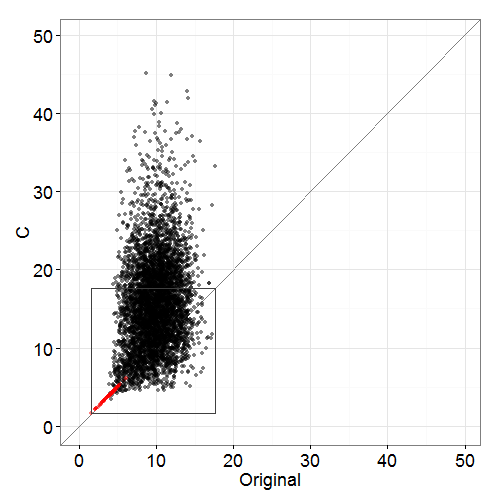}\\
\includegraphics[width=0.17\textwidth]{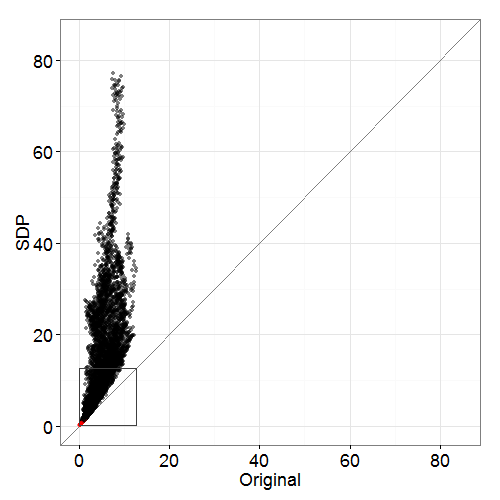} &
\includegraphics[width=0.17\textwidth]{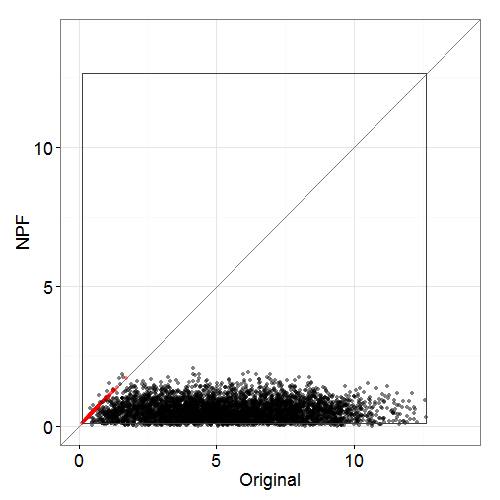}  &
\includegraphics[width=0.17\textwidth]{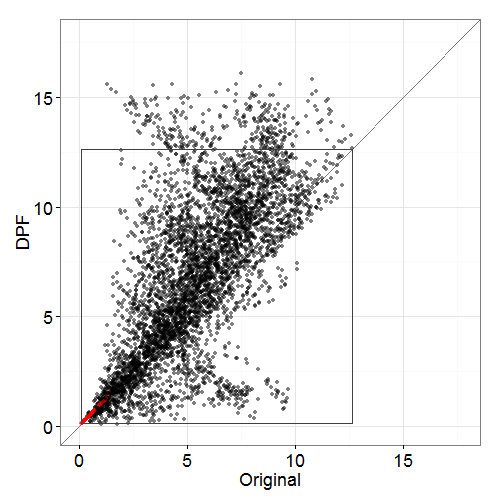}  & 
\includegraphics[width=0.17\textwidth]{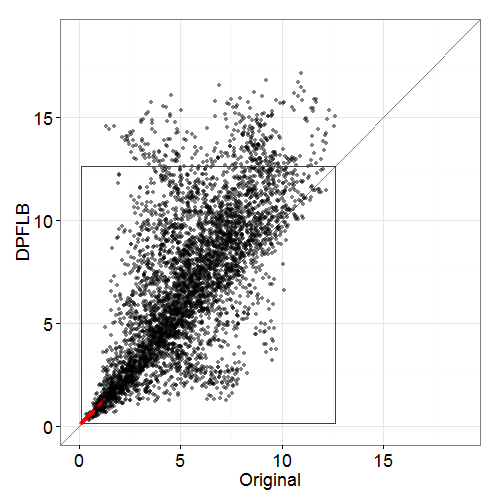} &
\includegraphics[width=0.17\textwidth]{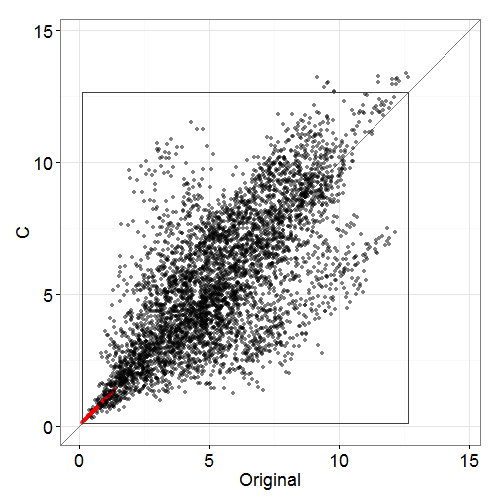}  \\
\footnotesize{(a) \texttt{SDP}} &
\footnotesize{(b) \texttt{NPF}} &
\footnotesize{(c) \texttt{DPF}} &
\footnotesize{(d) \texttt{DPFLB}} &
\footnotesize{(e) \texttt{C}} \\
\end{tabular}
\end{center}
\caption{\emph{Plots of $(d_{ij}, \widehat{d}_{ij})$ for all $i < j$ for a single reconstruction when the matrix to be completed contained only minimal spanning tree distances (shown in red);  the $y=x$ line indicates perfect matching; the box in each plot shows the extent of the original distances.  Results for all five methods are shown in each column; rows, from bottom to top, show increasing dimensionality of $p=2, 6, 10$. }}
\label{fig:UnifDistances}
\end{figure}

Figure \ref{fig:UnifDistances}
shows the completed distances $\widehat{d}_{ij}$ of each of the five methods paired with the actual distances $d_{ij}$ for one completion of matrix 1 having dimension $p= 2, 6, 10$.  The box within each plot shows the extent of the distances $d_{ij}$ and is identical in absolute magnitude across all plots.  The five methods appear as columns; the three rows show increasing dimensionality $p=2, 6, 10$ from bottom to top.

Consider the bottom row where $p=2$.  Here four of the five methods produce most, if not all, distances within the box given by the range of the original distances $d_{ij}$.   Of the two mst-preserving completions, \texttt{C} gives distances that are roughly symmetric about the $y=x$ line and clustered near it; a relatively few large distance appear outside the box at the top right.  In contrast, \texttt{DPFLB} produces distances that are more often above the $y=x$ line than below and are outside the box at top all along the range of the original distances.  The completion \texttt{DPF} is similar to \texttt{DPFLB} but will not necessarily reproduce the minimal spanning tree.

Most unusual in the bottom row are \texttt{SDP} and \texttt{NPF}.  Already for $p=2$, \texttt{SDP} produces huge distances  $\widehat{d}_{ij}$, far outside the box.  Moving up column (a) as the dimensionality increases, \texttt{SDP} produces larger and larger distances, much larger than the original distances and larger than those produced by any other completion method.  In contrast, \texttt{NPF} produces small distances, essentially all lying below the $y=x$ line.  This continues to be the case as the dimension increases;  distances produced by \texttt{NPF} ever smaller as $p$ increases with many becoming much smaller than any of the original distances $d_{ij}$.   This is the opposite of each of the other four methods whose distances become larger and larger as $p$ increases.  This contrast explains why in Figure \ref{fig:UnifPerformance} \texttt{NPF} showed little variation in $RDD$ compared to the other methods.   Because \texttt{NPF} produces small distances bounded below by 0,  $RDD$  is bounded for \texttt{NPF} but not for the others.

For all dimensions the constructive method \texttt{C} produces distances that are closer to the original distances than any of the other four methods. 

\subsection{Reproducing the minimal spanning tree}\label{sec:Uniform:repmst}
Both \texttt{C} and \texttt{DPFLB} reproduce the minimal spanning tree by design, as does \texttt{SDP}, except in this case largely by accident.  As Figure \ref{fig:UnifDistances} shows, the distances produced by \texttt{SDP} are typically so large that they do not change the minimal spanning tree.  The same effect can be seen with \texttt{DPF} where, because ever larger distances are produced as $p$ increases, the proportion of the minimal spanning tree retained by  \texttt{DPF} completions increases with $p$.  Completions by \texttt{NPF} on the other hand do not preserve the minimal spanning tree.  The small  distances produced by \texttt{NPF} interfere with the minimal spanning tree.

\subsection{Reproducing the point configurations}\label{sec:Uniform:Configurations}
To compare point configurations, we consider a completion from each of the five methods for matrix 1 when $p=2$ and $p=6$.  The data were first rotated to their principal coordinates as described in Section \ref{sec:Iris:PointConfigs}.

\begin{figure}[hbtp]
\begin{center}
\begin{tabular}{ccc}
\includegraphics[width=0.3\textwidth]{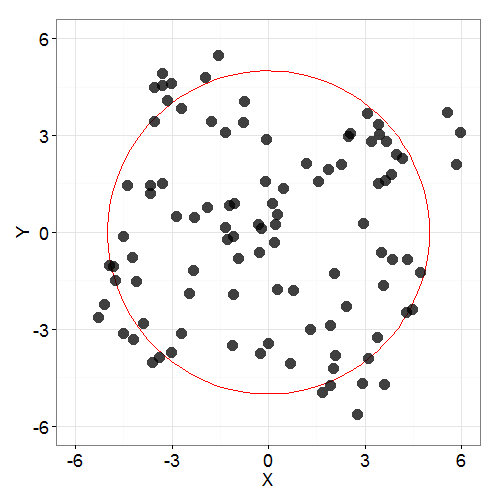} &
\includegraphics[width=0.3\textwidth]{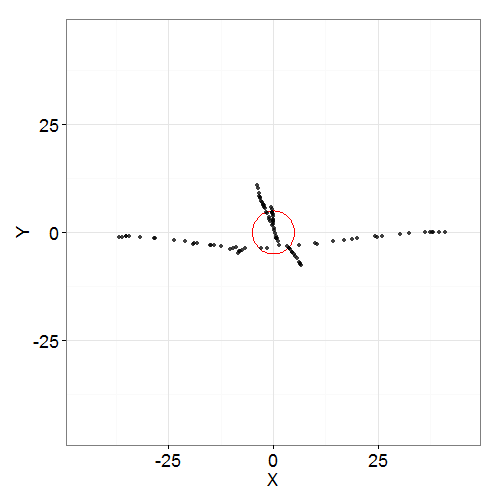} &
\includegraphics[width=0.3\textwidth]{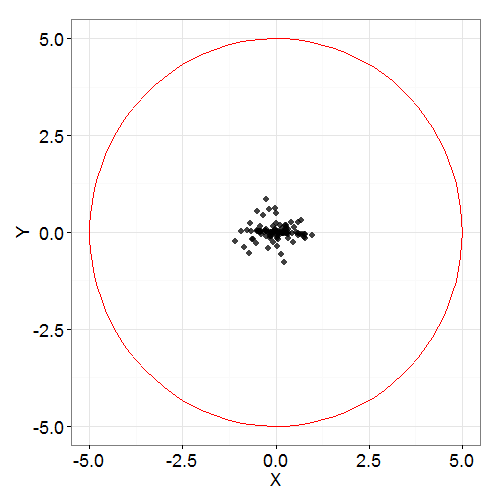} \\
\footnotesize{(a) Original}&
\footnotesize{(b) \texttt{SDP}} &
\footnotesize{(c) \texttt{NPF}} \\
\includegraphics[width=0.3\textwidth]{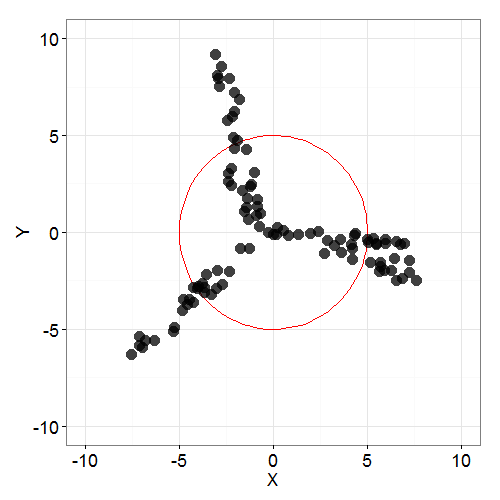} &
\includegraphics[width=0.3\textwidth]{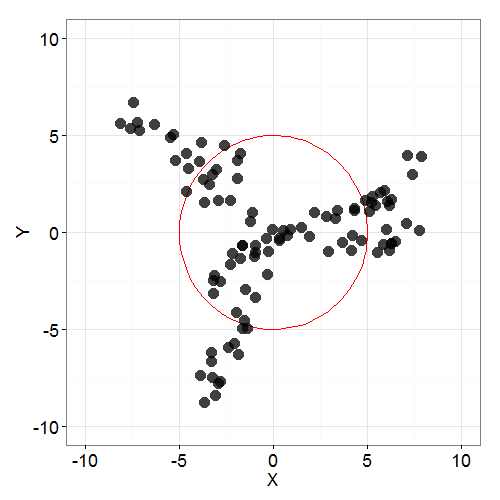} &
\includegraphics[width=0.3\textwidth]{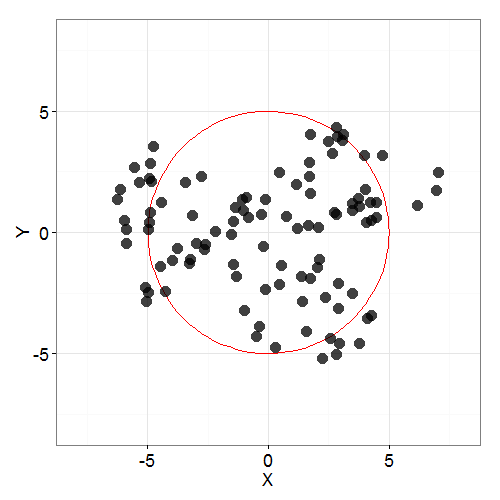} \\
\footnotesize{(d) \texttt{DPF}} &
\footnotesize{(e) \texttt{DPFLB}} &
\footnotesize{(f) \texttt{C}} \\
\end{tabular}
\end{center}
\caption{\emph{Uniform two dimensional data: the original and those reconstructed by the five methods. Coordinates are the principal coordinates for each configuration.  The red circle has the same centre and diameter in all plots. }}
\label{fig:UnifPointConfigs2D}
\end{figure}

Figure \ref{fig:UnifPointConfigs2D}
shows the original data in (a) and the reconstructed configurations in (b)--(f); each red circle has the same centre and radius in each plot.  Figure \ref{fig:UnifPointConfigs2D} reproduces the findings from the bottom row of Figure \ref{fig:UnifDistances} in that \texttt{SDP} produces unusually large distances, \texttt{NPF} unusually small distances, and the mst-preserving methods \texttt{DPFLB} and \texttt{C} produce distances closer to those of the original data, with \texttt{C} being the closest.

As was the case with the Iris data, Figure \ref{fig:UnifPointConfigs2D} shows again that the completion methods tend to concentrate points near lines.  The mst-preserving \texttt{DPFLB} spreads the configuration out more than does \texttt{DFP} but not nearly as much as does \texttt{C}.  For $p=2$ \texttt{C} produces a point configuration that is more like the original data than any of the others.

\begin{figure}[hbtp]
\begin{tabular}{cc}
\includegraphics[width=0.5\textwidth]{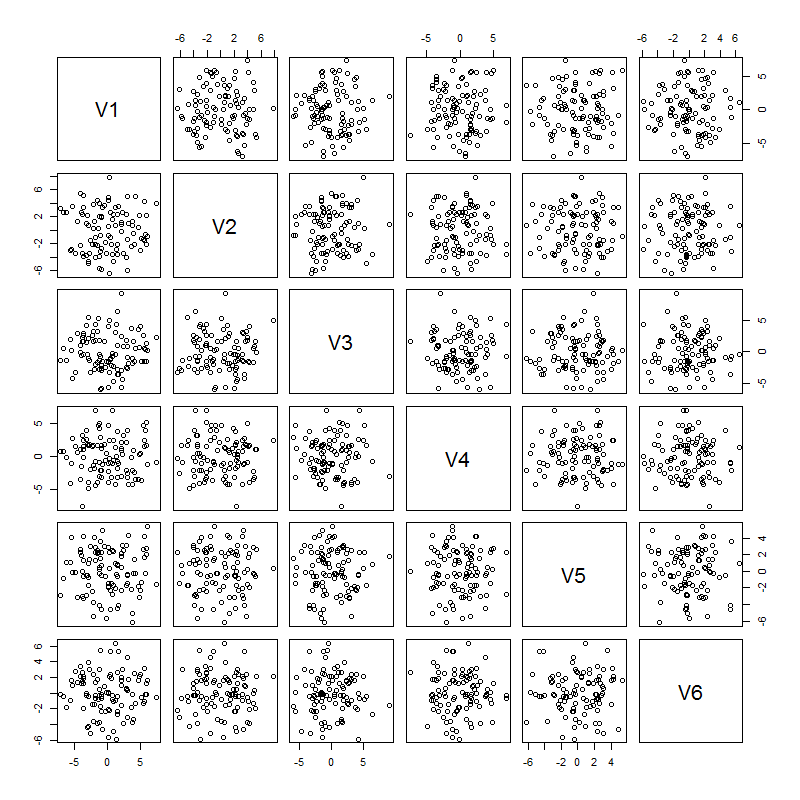} &
\includegraphics[width=0.5\textwidth]{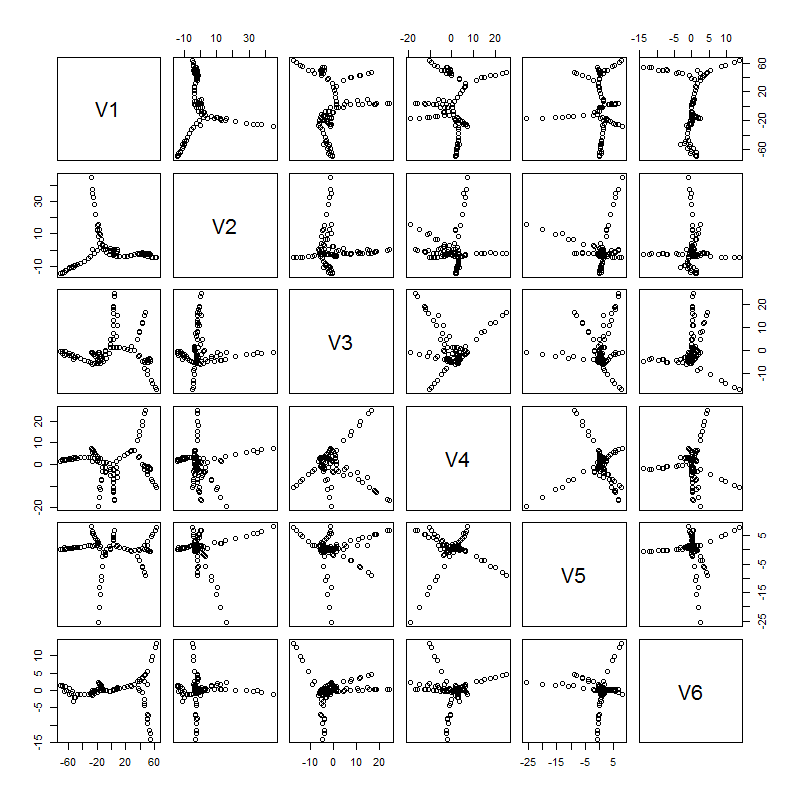} \\
\footnotesize{(a) Original}&
\footnotesize{(b) \texttt{SDP}} \\
\includegraphics[width=0.5\textwidth]{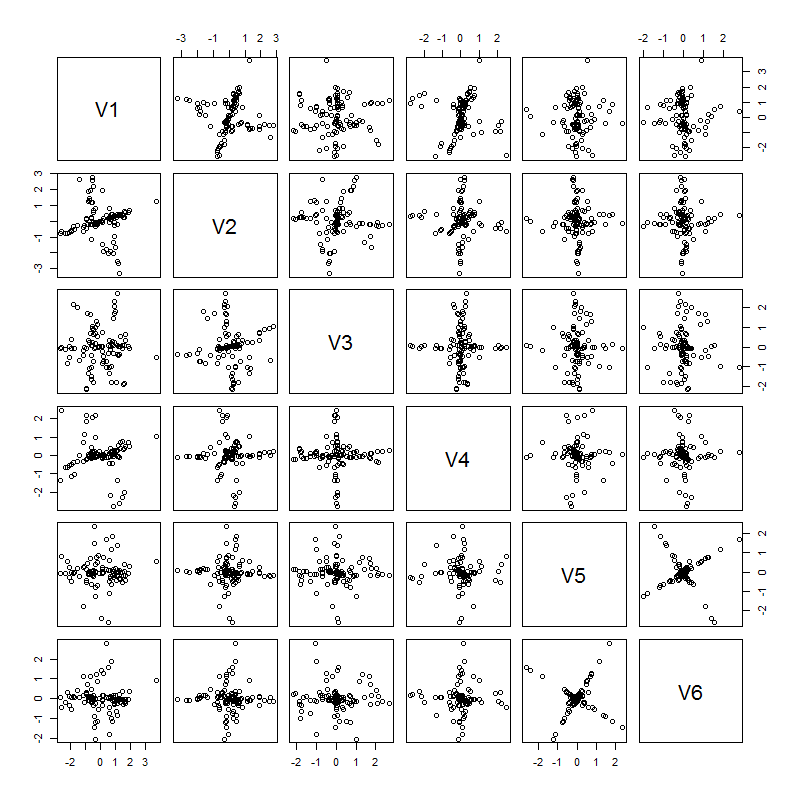} &
\includegraphics[width=0.5\textwidth]{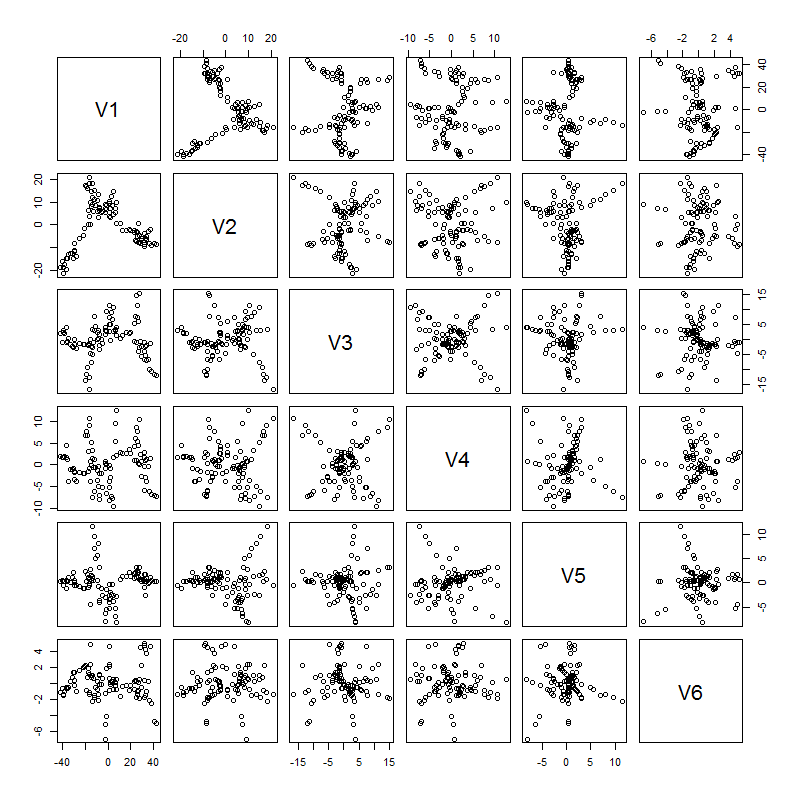} \\
\footnotesize{(c) \texttt{NPF}}&
\footnotesize{(d) \texttt{DPF}} \\
\includegraphics[width=0.5\textwidth]{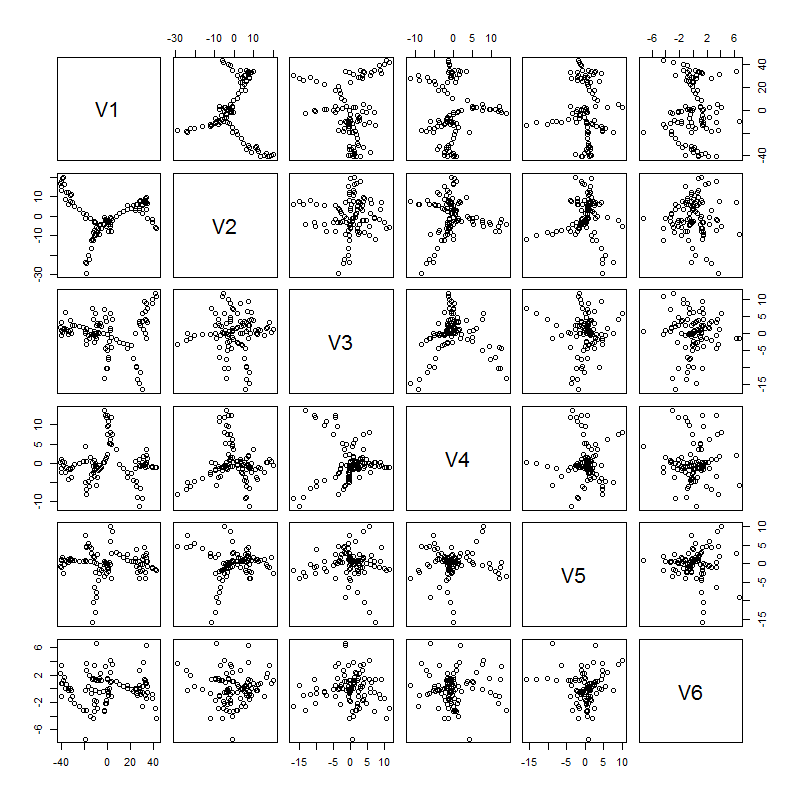} &
\includegraphics[width=0.5\textwidth]{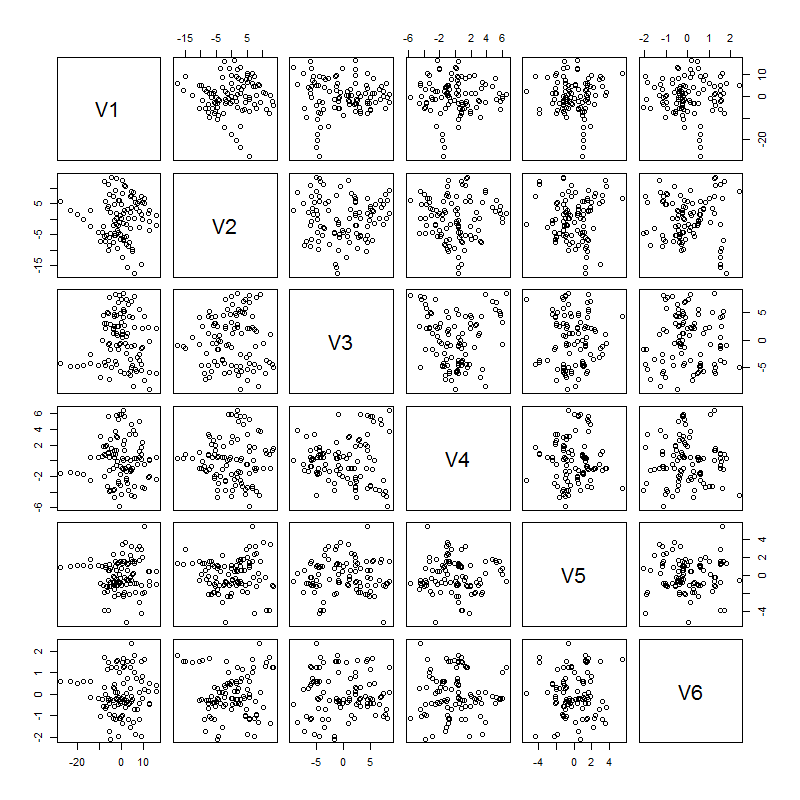} \\
\footnotesize{(e) \texttt{DPFLB}}&
\footnotesize{(f) \texttt{C}} \\
\end{tabular}\caption{\emph{Uniform within a six-dimensional hypercube. Data point configurations reconstructed in six dimensions using each reconstruction's (including the original data) principal coordinates. }}
\label{fig:Unif6DScatters}
\end{figure}

Figure \ref{fig:Unif6DScatters}
shows the point configurations produced when $p=6$.  Much the same patterns prevail as were seen earlier with $p=2$ and also when $p=4$ for the Iris data.  \texttt{SDP} has large distances and strongly linear configurations;  \texttt{NPF} has small distances but also exhibits star-shaped linear structure; \texttt{DPF} and \texttt{DPFLB} also show linear structure and some star shape;  \texttt{C} shows a dispersed configuration most like the original data but also has outlying points in the $V1$ direction.

\section{Concluding remarks}\label{sec:Conclusions}
Interest here has been to construct point configurations in some embedding dimension $p$ given the matrix of (squared) Euclidean distances corresponding to the minimal spanning tree.  All other distances are missing  and we would like to have the minimal spanning tree preserved in any configuration produced.   As outlined in Sections \ref{sec:edm} and \ref{sec:mst}, this is a special case of the Euclidean distance matrix completion problem.  

Section \ref{sec:mst} showed how the problem could be cast in the formulation of \cite{Trosset} and in Section \ref{sec:dbflb} provided Algorithm \ref{alg:mstLB} to implement the solution.  Section \ref{sec:con} took a more constructive approach.  Here no optimization problem or solution was used.  Instead a solution was constructed through a guided random search, as implemented in Algorithm \ref{alg:constructive}.   In this way a multiplicity of solutions could be produced.

In Section \ref{sec:expt}, all methods were compared on reconstructions of known point configurations.  In both experimental setups, reproducing the Iris data in Section \ref{sec:Iris} or the uniform data in Section \ref{sec:uniform},  the constructive method \texttt{C} was observed to outperform all others.  It preserved the minimal spanning tree, was orders of magnitude faster than all others,  more accurately reproduced the original (non-mst) inter-point distances, and, perhaps most importantly, gave point configurations much more like the original data.  Note also that all other methods (including the mst-preserving \texttt{DPFLB})  produced unnaturally regular shaped point configurations as evidenced by Figures \ref{fig:IrisScatters}, \ref{fig:UnifPointConfigs2D}, and \ref{fig:Unif6DScatters}; any predisposition to do so is a concern. 

In contrast to minimal spanning tree completion, when 75\% or fewer of the distances were missing, all three Euclidean distance matrix completion methods (viz. \texttt{SDP}, \texttt{NPF}, and \texttt{DPF}) performed fairly well in terms of accuracy (see bottom plot of Figure \ref{fig:IrisPerformance}) performed well in terms of accuracy as measured by $RDD$  ($RDD < 0.01$ for the most part for \texttt{SDP} but $RDD < 0.000001$ for \texttt{DPF} and $RDD < 0.0000001$ for \texttt{NPF}).  Beyond 75 or 80\% missing the accuracy of the three methods degraded quickly.  The more significant contrast between the three was the time taken to completion (see top plot of Figure \ref{fig:IrisPerformance}).  Here \texttt{SDP} took 1-5 orders of magnitude longer than did either \texttt{DPF} or \texttt{NPF} for all percentages of missing.  With the exception of the minimal spanning tree case, \texttt{DPF} was fastest (from 10 to about 100 times faster than \texttt{NPF}).

The case of the minimal spanning tree distance matrix completion is unusual in at least three ways.  
First, the completed matrix must satisfy an additional constraint on its structure -- it must preserve the minimal spanning tree.  Second, nearly all of the distances are missing which enlarges the solution space and is more computationally challenging for both \texttt{NPF} and \texttt{DPF} (\texttt{SDP} is about the same but is also far more computationally intensive than either of \texttt{NPF} or \texttt{DPF} in all cases). Finally, the initial matrix forms only a spanning tree and so forces little structure on the positioning of points compared to a denser graph.   The last two expand the solution space, while the first constrains it.   The expansion seems to dominate the constraint here in that there are so many constraint satisfying solutions that a guided search like that given by the constructive method \texttt{C} has no trouble finding a viable one.  And it does so very quickly, allowing a multiplicity of solutions to potentially be examined and compared.

As Section \ref{sec:edm} reviewed, the Euclidean distance matrix completion problem is usually treated as a constrained optimization problem.  But why?  In the case of minimal spanning tree distance matrix completion, the  remarkable out-performance of the other methods by \texttt{C} suggests that this might not always be the best approach.  Perhaps this observation applies more broadly.   After all, none of the methods described in Section \ref{sec:edm} even find a global optimum; they all find a local optimum  (the semi-definite programming formulation is convex only when the embedding dimension $p$ is unrestricted and allowed to be as large as $n$).  In principle, then, multiple local optima could be found for each method which is good since, for the minimal spanning tree problem at least, there are in fact many global optima to choose from.

Instead of focusing on completing the Euclidean distance matrix, we could focus on producing a point configuration as does \texttt{NPF} of \cite{OLeary} who arrive at the minimization problem given by Equation \eqref{eq:pcf} from Section \ref{sec:npf}.  The objective function they minimize, from Equation \ref{eq:pcf:objectivefn}, is
\[
f_{A, D}(\m{X}) = \sum_{i=1}^n \sum_{j=1}^n a_{ij}^2 (d_{ij} - \norm{\ve{x}_i - \ve{x}_j}^2)^2.
\]
A constructive approach would be to simply choose positions $\widehat{\ve{x}}_i$ such that $d_{ij} = \norm{\widehat{\ve{x}}_i - \widehat{\ve{x}}_j}^2$ whenever $a_{ij}=1$; the value of $ \norm{\widehat{\ve{x}}_i - \widehat{\ve{x}}_j}$ doesn't matter whenever $a_{ij} = 0$.  There is typically a multiplicity of solutions $\widehat{\m{X}}$ which will yield  $f_{A, D}(\widehat{\m{X}}) = 0$, we need only find one.   A random search might do. There are no doubt numerous problem instances where particular graph structure given by $\m{A}$ can be exploited to quickly find solutions through random search; or a combination of random search and, say, \texttt{NPF}, for the denser parts of the graph.

Finally, when there are multiple optima that can be easily found, we might reasonably choose between them on other grounds.  For example, we have already begun work on adaptations to Algorithm \ref{alg:constructive} that, for data visualization applications, restrict the shape of the minimal spanning tree.
Given that \texttt{C} can find  $100$s or even tens of thousands of solutions in the time taken by the other methods to find one, some further constraint on the random positions proposed is easily practicable.  There are no doubt many interesting constraints on the shape of configurations for the minimal spanning tree completion problem as well as graph structures other than the minimal spanning tree which might  be exploited to guide a random search. 

\appendix

\section{Proof of the convex cone}

\theorem{The set
${\cal M}_n := {\cal M}_n ( \m{A}, \m{A}^{\star}, \sm{\Delta}^{\star} ) = \{ \sm{\Delta} \in {\cal C}_n : amst(\m{A}, \sm{\Delta}) =   amst(\m{A}^{\star}, \sm{\Delta}^{\star}) \}
$ is a convex cone.
}
\proof{Two things need to be shown: 
\begin{enumerate}
\item if $\sm{\Delta} \in {\cal M}_n$ then $\alpha \sm{\Delta} \in {\cal M}_n$ for any real $\alpha > 0$, and
\item if $\sm{\Delta} \in {\cal M}_n$ and $\sm{\Lambda} \in {\cal M}_n$ then 
$\sm{\Gamma} = \alpha \sm{\Delta} + \beta \sm{\Lambda} \in {\cal M}_n$ for any reals $\alpha, \beta > 0$.
\end{enumerate}
First, note that ${\cal C}_n$ is clearly a convex cone from its definition, so the requirement that members of ${\cal M}_n$ also be members of ${\cal C}_n$ is trivially satisfied for both items above, so we need only check that the minimal spanning tree requirements are met.

Item 1 is also trivially true.  If $amst(\m{A}, \sm{\Delta}) =   amst(\m{A}^{\star}, \sm{\Delta}^{\star})$ then a common rescaling of all elements in $\sm{\Delta}$ will make no change to the adjacency
of any minimal spanning tree.

Item 2 is proved by showing that $amst(\m{A},\sm{\Delta}) \subseteq amst(\m{A},\sm{\Gamma})$, then that  \\ $amst(\m{A},\sm{\Gamma}) \subseteq amst(\m{A},\sm{\Delta})$, implying $amst(\m{A},\sm{\Gamma}) = amst(\m{A},\sm{\Delta})$.

To show $amst(\m{A}, \sm{\Delta}) \subseteq amst(\m{A},\sm{\Gamma})$, let $\m{M} \in amst(\m{A}, \sm{\Delta})$ and $\m{B} \in {\cal A}_n^{\star}$ be any spanning tree of both $\m{A}$ and $\m{A}^{\star}$.  We write (twice) the sum of dissimilarities of $\m{M}\Had \sm{\Gamma}$ as
\begin{align*}
\tr{\ve{1}}(\m{M} \Had \sm{\Gamma}) \ve{1}  &
=  \tr{\ve{1}}(\m{M} \Had (\alpha \sm{\Delta})) \ve{1} + \tr{\ve{1}}(\m{M} \Had (\beta \sm{\Lambda})) \ve{1} \\
&
\le \tr{\ve{1}}(\m{B} \Had (\alpha \sm{\Delta})) \ve{1} + \tr{\ve{1}}(\m{B} \Had (\beta \sm{\Lambda})) \ve{1} \\
& = \tr{\ve{1}}(\m{B} \Had \sm{\Gamma}) \ve{1}
\end{align*}
which implies that $\m{M} \in amst(\m{A},\sm{\Gamma})$ and so $amst(\m{A},\sm{\Delta}) \subseteq amst(\m{A},\sm{\Gamma})$.

To show $amst(\m{A},\sm{\Gamma}) \subseteq amst(\m{A},\sm{\Delta})$, let $\m{M}\in amst(\m{A},\sm{\Delta})$ and $\m{B} \in amst(\m{A}, \sm{\Gamma})$.  We prove the result by contradiction.  
Suppose $\m{B} \not\in amst(\m{A}, \sm{\Delta})$.  
Then we have both that  
$\tr{\ve{1}}(\m{B} \Had (\alpha \sm{\Delta})) \ve{1} >  \tr{\ve{1}}(\m{M} \Had (\alpha \sm{\Delta})) \ve{1}$ 
and that 
$\tr{\ve{1}}(\m{B} \Had (\beta \sm{\Lambda})) \ve{1} > \tr{\ve{1}}(\m{M} \Had (\beta \sm{\Lambda})) \ve{1}$.  
Together these imply that
\begin{align*}
\tr{\ve{1}}(\m{B} \Had \sm{\Gamma}) \ve{1}  &
=  \tr{\ve{1}}(\m{B} \Had (\alpha \sm{\Delta})) \ve{1} + \tr{\ve{1}}(\m{B} \Had (\beta \sm{\Lambda}) )\ve{1} \\
&
>  \tr{\ve{1}}(\m{M} \Had (\alpha \sm{\Delta})) \ve{1} + \tr{\ve{1}}(\m{M} \Had (\beta \sm{\Lambda}) )\ve{1} \\
& = \tr{\ve{1}}(\m{M} \Had \sm{\Gamma}) \ve{1}
\end{align*}
which means that $\m{M}$ has a shorter spanning tree and hence $\m{B} \not\in amst(\m{A}, \sm{\Gamma})$, a contradiction.  Therefore $\m{B} \in amst(\m{A}, \sm{\Delta})$ and hence
$amst(\m{A},\sm{\Gamma}) \subseteq amst(\m{A},\sm{\Delta})$.

}

\bibliographystyle{latexStyles/siamplain}
\bibliography{TheMSTProblem}

\end{document}